\documentclass[11pt]{article}

\usepackage{caption2}

\usepackage{graphicx}

\usepackage[numbers]{natbib}

\usepackage[all,cmtip]{xy}

\usepackage{setspace}
\usepackage{mathtools, bbm,youngtab}
\usepackage{amscd,mathrsfs}
\usepackage{amssymb,amsmath,dsfont,amsfonts,amsthm}
\usepackage{euscript,enumerate}

\usepackage{stmaryrd}

\usepackage[cbgreek]{textgreek} 

\usepackage{booktabs}
\renewcommand{\arraystretch}{1.15} 
\usepackage{multirow}

\usepackage{hhline} 

\usepackage{cases} 

\usepackage{arydshln} 

\usepackage{color}

\DeclareMathAlphabet{\mathpzc}{OT1}{pzc}{m}{it}

\usepackage{epstopdf}

\allowdisplaybreaks 

\usepackage{mathtools}                      
\mathtoolsset{showonlyrefs=true}          

 \theoremstyle{plain}
 \newtheorem{thm}{Theorem}
 \newtheorem{cor}[thm]{Corollary}

 \theoremstyle{definition}

 \numberwithin{equation}{section}

\usepackage{hyperref}
\hypersetup{colorlinks=true, linkcolor=blue}
\hypersetup{colorlinks=true,citecolor=blue}

\makeatletter
\DeclareFontFamily{OMX}{MnSymbolE}{}
\DeclareSymbolFont{MnLargeSymbols}{OMX}{MnSymbolE}{m}{n}
\SetSymbolFont{MnLargeSymbols}{bold}{OMX}{MnSymbolE}{b}{n}
\DeclareFontShape{OMX}{MnSymbolE}{m}{n}{
    <-6>  MnSymbolE5
   <6-7>  MnSymbolE6
   <7-8>  MnSymbolE7
   <8-9>  MnSymbolE8
   <9-10> MnSymbolE9
  <10-12> MnSymbolE10
  <12->   MnSymbolE12
}{}
\DeclareFontShape{OMX}{MnSymbolE}{b}{n}{
    <-6>  MnSymbolE-Bold5
   <6-7>  MnSymbolE-Bold6
   <7-8>  MnSymbolE-Bold7
   <8-9>  MnSymbolE-Bold8
   <9-10> MnSymbolE-Bold9
  <10-12> MnSymbolE-Bold10
  <12->   MnSymbolE-Bold12
}{}

\let\llangle\@undefined
\let\rrangle\@undefined
\DeclareMathDelimiter{\llangle}{\mathopen}%
                     {MnLargeSymbols}{'164}{MnLargeSymbols}{'164}
\DeclareMathDelimiter{\rrangle}{\mathclose}%
                     {MnLargeSymbols}{'171}{MnLargeSymbols}{'171}
\makeatother

\title{\Large{\textbf{Four-Field Mixed Finite Element Methods for\\ Incompressible Nonlinear Elasticity}}}

\author{Arzhang Angoshtari\thanks{Department of Civil and Environmental Engineering, The George Washington University, Washington, DC 20052. E-mail: aangoshtari@gwu.edu.}}

\begin{document}

\maketitle

\begin{abstract} We introduce conformal mixed finite element methods for $2$D and $3$D incompressible nonlinear elasticity in terms of displacement, displacement gradient, the first Piola-Kirchhoff stress tensor, and pressure, where finite elements for the $\mathrm{curl}$ and the $\mathrm{div}$ operators are used to discretize strain and stress, respectively. These choices of elements follow from the strain compatibility and the momentum balance law. Some inf-sup conditions are derived to study the stability of methods. By considering $96$ choices of simplicial finite elements of degree less than or equal to $2$ in $2$D and $3$D, we conclude that $28$ choices in $2$D and $6$ choices in $3$D satisfy these inf-sup conditions. The performance of stable finite element choices are numerically studied. Although the proposed methods are computationally more expensive than the standard two-field methods for incompressible elasticity, they are potentially useful for accurate approximations of strain and stress as they are independently computed in the solution process.
\end{abstract}

\begin{description}
\item[Keywords.] Incompressible nonlinear elasticity; mixed finite element methods; inf-sup conditions.
\end{description}


\section{Introduction}

It is well-known that incompressible nonlinear elasticity is useful to describe the mechanical behavior of various soft tissues \citep{Hu2003}. Deriving stable numerical methods for modeling such tissues that may undergo large deformations is a challenging task, for example see \citep{AuDaLoReTaWr2013} and references therein. In this paper, we develop four-field mixed finite element methods for incompressible nonlinear elasticity by extending three-field mixed finite element methods introduced in \citep{AnGr2019I} for compressible nonlinear elasticity. These four-field methods are based on a mixed formulation with displacement, displacement gradient, the first Piola-Kirchhoff stress tensor, and pressure as the independent unknowns. A potential advantage of this mixed formulation is providing accurate approximations of strain and stress, since they are independent unknowns and are explicitly calculated in the solution process, that is, no post-processing is required to indirectly compute strain and stress from displacement.       

Similar to \citep{AnGr2019I}, finite elements suitable for the $\mathrm{curl}$ and the $\mathrm{div}$ operators are respectively employed to discretize displacement gradient and stress in the proposed conformal mixed finite element methods. One can show that these choices follow from the strain compatibility and the balance law \citep{AngoshtariYavari2014I,AngoshtariYavari2014II}. To discretize displacement and pressure, $H^{1}$- and $L^{2}$-conformal finite element spaces are used, respectively.

A stability analysis based on the abstract theory of \citep{PoRa1994,CaRa1997} for Galerkin approximations of nonlinear problems is provided. More specifically, an inf-sup condition is written which is locally sufficient for the convergence of solutions of nonlinear finite element methods to a regular solution of incompressible nonlinear elasticity. This inf-sup condition can be interpreted as the necessary and sufficient condition for the well-posedness of the linear problems associated to Newtons' iterations for solving nonlinear finite element methods. Based on this interpretation, we also write $6$ other inf-sup conditions that are necessary for the stability of Newtons' iterations. Using these conditions, we derive some relations between the dimension of finite element spaces for different unknowns, which are necessary for the stability of Newtons' iterations. 

By considering $96$ choices of simplicial finite elements of degree less than or equal to $2$ for obtaining mixed finite element methods, we conclude that $68$ choices in $2$D and $90$ choices in $3$D do not satisfy all these inf-sup conditions, in general, and lead to unstable methods. The performance of stable choices in $2$D and $3$D are studied by solving numerical examples.

This paper is organized as follows: The four-field mixed formulation and the associated conformal finite element methods are discussed in Sections 2 and 3. Newtons' iterations for solving nonlinear finite element methods are provided in Section 4. In Section 5, different inf-sup conditions are derived for studying the stability of finite element methods and it is shown that some choices of finite elements lead to unstable methods as they violate the inf-sup conditions. In Section 6, the inf-sup conditions are numerically studied. Moreover, by solving two numerical examples, the performance of stable finite element methods are investigated in $2$D and $3$D. Finally, some concluding remarks are given in Section 7.

\section{A Mixed Formulation for Incompressible Elasticity}

We write a four-field mixed formulation for incompressible nonlinear elasticity in this section. For simplicity, we consider static problems here, however, this formulation can be readily extended to time-dependent problems as well. Let $\mathcal{B}$ denote the reference configuration of a $2$D or $3$D elastic body and let $\boldsymbol{N}$ be the unit outward normal vector field associated to the boundary $\partial\mathcal{B}=\Gamma_{1}\cup\Gamma_{2}$ of $\mathcal{B}$.
 
Suppose $(\boldsymbol{U},\boldsymbol{K}, \boldsymbol{P})$ are respectively displacement, displacement gradient, and the first Piola-Kirchhoff stress tensor. The constitutive relation of incompressible elastic bodies can be expressed as  $\boldsymbol{P} = \mathbb{P}(\boldsymbol{K}) - p\cdot(\boldsymbol{I}+\boldsymbol{K})^{-T}$, where $\mathbb{P}$ is a given function of $\boldsymbol{K}$, the real-valued function $p$ is pressure, and $\boldsymbol{I}$ is the identity tensor \citep[Section $30$]{TrNo65}. The boundary value problem of incompressible nonlinear elastostatics can be stated as:

\bigskip
\begin{minipage}{.95\textwidth}{\textit{Given a body force $\boldsymbol{B}$, a displacement $\overline{\boldsymbol{U}}$ of $\mathcal{B}$, and a traction vector field $\overline{\boldsymbol{T}}$ on $\Gamma_{2}$, find $(\boldsymbol{U},\boldsymbol{K},\boldsymbol{P},p)$ such that}}
\begin{align}
\left. \begin{array}{l}
\mathbf{div}\,\boldsymbol{P} = - \boldsymbol{B}, \\
\boldsymbol{K}-\mathbf{grad}\,\boldsymbol{U}=\boldsymbol{0}, \\
\boldsymbol{P}-\mathbb{P}(\boldsymbol{K})+ p(\boldsymbol{I}+\boldsymbol{K})^{-T}=\boldsymbol{0}, \\
\det (\boldsymbol{I}+\boldsymbol{K})=1,
\end{array} \right\}  \text{ in }\mathcal{B}, \label{Strong} \\
\boldsymbol{U}=\overline{\boldsymbol{U}},  \text{ on } \Gamma_{1}, \quad 
\boldsymbol{P}(\boldsymbol{N})=\overline{\boldsymbol{T}},  \text{ on } \Gamma_{2}. \label{BC}
\end{align}
\end{minipage}
\bigskip

To write a week form based on the above equations, we consider the following spaces: Let $L^{2}(\mathcal{B})$ be the standard space of square-integrable functions and let $H^{1}(\mathcal{B})$ be the Sobolev space of $L^{2}$-functions with first-order derivatives of $L^{2}$-class. The space of vector fields in $\mathbb{R}^{n}$, $n=2,3$, with components of $H^{1}$-class is denoted by $[H^{1}(\mathcal{B})]^{n}$. The space of $H^{1}$-vector fields that vanish on $\Gamma_{1}$ is written as $[H^{1}_{1}(\mathcal{B})]^{n}$. 

We also consider the spaces $H^{\mathbf{c}}(\mathcal{B})$ and $H^{\mathbf{d}}(\mathcal{B})$, which are the spaces of second-order tensor fields with $L^{2}$-components that their $\mathbf{curl}$ and $\mathbf{div}$ are respectively of $L^{2}$-class. Consider the $n\times n$ matrix $[\boldsymbol{S}]$ representing the Cartesian components of a second-order tensor field $\boldsymbol{S}$ in $\mathbb{R}^{n}$. Since $\mathbf{curl}\,\boldsymbol{S}$ and $\mathbf{div}\,\boldsymbol{S}$ are respectively obtained by applying the standard $\mathrm{curl}$ and $\mathrm{div}$ of vector fields to each row of $[\boldsymbol{S}]$, the spaces $H^{\mathbf{c}}(\mathcal{B})$ and $H^{\mathbf{d}}(\mathcal{B})$ can be considered as $n$-copies of the corresponding spaces for vector fields.    

For obtaining a concise form for the weak formulation, we employ the following notation: The standard inner product of $\mathbb{R}^{n}$ is denoted by ``$\boldsymbol{\cdot}$'' and the $L^{2}$-inner products of real-valued functions, vector fields, and tensor fields are denoted by $\llangle,\rrangle$, that is, $\llangle f,g\rrangle:=\int_{\mathcal{B}}fg\,dV$, $\llangle\boldsymbol{Y},\boldsymbol{Z}\rrangle:=\int_{\mathcal{B}}Y^{I}Z^{I}dV$, and $\llangle\boldsymbol{S},\boldsymbol{T}\rrangle:=\int_{\mathcal{B}}S^{IJ}T^{IJ}dV$, where we consider the summation convention on repeated indices. By taking the $L^{2}$-inner product of the equations \eqref{Strong} with suitable test functions and using Green's formula, one can write the following four-field mixed formulation for incompressible nonlinear elasticity:

\bigskip
\begin{minipage}{.95\textwidth}{\textit{Given a body force $\boldsymbol{B}$, a boundary displacement $\overline{\boldsymbol{U}}$, and a surface traction vector field $\overline{\boldsymbol{T}}$ on $\Gamma_{2}$, find $(\boldsymbol{U},\boldsymbol{K},\boldsymbol{P},p)\in [H^{1}(\mathcal{B})]^{n}\times H^{\mathbf{c}}(\mathcal{B}) \times H^{\mathbf{d}}(\mathcal{B})\times L^{2}(\mathcal{B})$ such that $\boldsymbol{U}= \overline{\boldsymbol{U}}$, on $\Gamma_1$, and }}
\begin{equation}\label{IncompWeakForm}
\begin{alignedat}{3}
&\llangle\boldsymbol{P},\mathbf{grad}\,\boldsymbol{\Upsilon}\rrangle =\llangle \boldsymbol{B},\boldsymbol{\Upsilon} \rrangle + \int_{\Gamma_{2}}\overline{\boldsymbol{T}}\boldsymbol{\cdot}\boldsymbol{\Upsilon} dA, &\quad &\forall \boldsymbol{\Upsilon}\in [H^{1}_{1}(\mathcal{B})]^{n},\\
&\llangle\mathbf{grad}\,\boldsymbol{U},\boldsymbol{\lambda}\rrangle - \llangle \boldsymbol{K},\boldsymbol{\lambda}\rrangle = 0,& &\forall \boldsymbol{\lambda}\in H^{\mathbf{c}}(\mathcal{B}),\\
&\llangle \mathbb{P}(\boldsymbol{K}),\boldsymbol{\pi}\rrangle - \llangle \boldsymbol{P}, \boldsymbol{\pi} \rrangle - \llangle p\,(\boldsymbol{I} +\boldsymbol{K})^{-T}, \boldsymbol{\pi} \rrangle  = 0, & &\forall \boldsymbol{\pi}\in H^{\mathbf{d}}(\mathcal{B}), \\
&\llangle \det(\boldsymbol{I} + \boldsymbol{K}), q \rrangle = \llangle 1,q\rrangle, & &\forall q\in L^{2}(\mathcal{B}).  
\end{alignedat}
\end{equation}
\end{minipage}
\bigskip

The above weak form is an extension of the formulation of \citep{AnGr2019I} for compressible elasticity to incompressible elasticity. The choices of the spaces $H^{\mathbf{c}}(\mathcal{B})$ and $H^{\mathbf{d}}(\mathcal{B})$ guarantee that physically relevant jump conditions, namely, the Hadamard jump condition for strain and the continuity of traction vector fields at interfaces for stress, will be hold in the weak sense and also on the discrete-level. Following the discussion of \citep[Remark 2]{AnGr2019I}, it is not hard to show that even for hyperelastic materials, the above weak formulation is not associated to a stationary point of any functional.

\section{Conformal Finite Element Methods}
Suppose $[V^{1}_{h}]^{n}$, $V^{\mathbf{c}}_{h}$, $V^{\mathbf{d}}_{h}$, and $V^{\mathrm{D}}_{h}$ are finite element spaces such that $[V^{1}_{h}]^{n}\subset [H^{1}(\mathcal{B})]^{n} $, $V^{\mathbf{c}}_{h} \subset H^{\mathbf{c}}(\mathcal{B})$, $V^{\mathbf{d}}_{h} \subset H^{\mathbf{d}}(\mathcal{B})$, and $V^{\mathrm{D}}_{h} \subset L^{2}(\mathcal{B})$. Also let $[V^{1}_{h,1}]^{n}= [V^{1}_{h}]^{n}\cap [H^{1}_{1}(\mathcal{B})]^{n}$ and let $\mathcal{I}^{1}_{h}$ be the interpolation operator associated to $[V^{1}_{h}]^{n}$. Then, one can write the following conformal mixed finite element method for \eqref{IncompWeakForm}:

\begin{figure}[t]
\begin{center}
\includegraphics[scale=.68,angle=0]{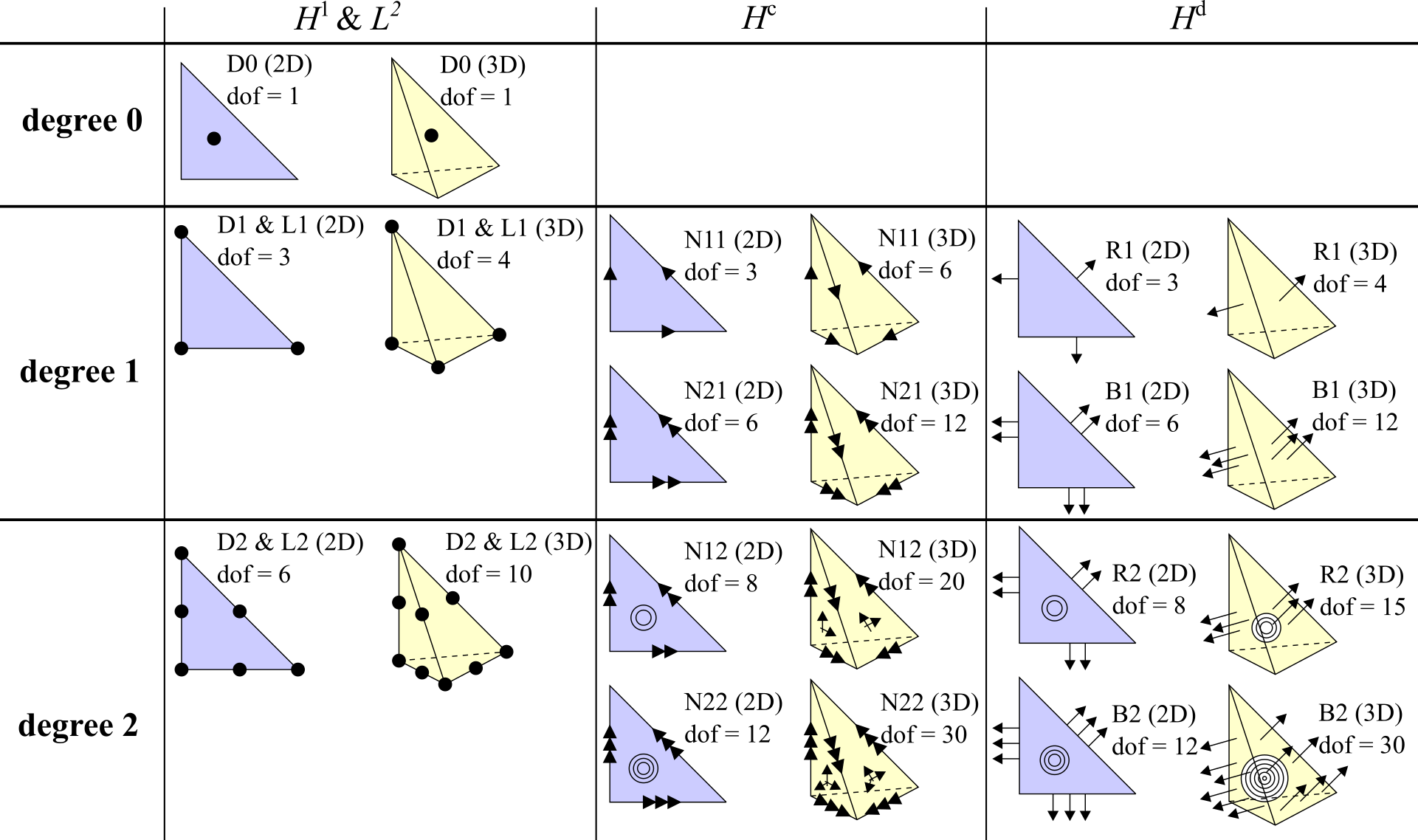}
\end{center}
\vspace*{-0.2in}
\caption{\footnotesize Conventional finite element diagrams for elements considered in this work: The discontinuous $L^{2}$ elements $\mathrm{D}i$ and the continuous $H^{1}$ elements $\mathrm{L}i$ correspond to the standard simplicial Lagrange elements of degree $i$. The $H^{c}$ element $\mathrm{N}{ji}$ stands for the $i$-th degree N\'{e}d\'{e}lec element of the $j$-th kind. The $H^{d}$ elements $\mathrm{R}i$ and $\mathrm{B}i$ respectively stand for the Raviart-Thomas element of degree $i$ and the Brezzi-Douglas-Marini element of degree $i$. Arrows parallel (normal) to an edge or a face denote degrees of freedom associated to tangent (normal) components of vector fields along that edge or face. Only degrees of freedom associated to visible edges and faces are shown.} 
\label{FEO012}
\end{figure}

\bigskip
\begin{minipage}{.95\textwidth}{\textit{Given a body force $\boldsymbol{B}$, a boundary displacement $\overline{\boldsymbol{U}}$, and a surface traction vector field $\overline{\boldsymbol{T}}$ on $\Gamma_{2}$, find $(\boldsymbol{U}_{h},\boldsymbol{K}_{h},\boldsymbol{P}_{h},p_{h})\in [V^{1}_{h}]^{n}\times V^{\mathbf{c}}_{h} \times V^{\mathbf{d}}_{h}\times V^{\mathrm{D}}_{h}$ such that $\boldsymbol{U}_{h}= \mathcal{I}^{1}_{h}(\overline{\boldsymbol{U}})$, on $\Gamma_1$, and }}
\begin{equation}\label{IncompFEMs}
\begin{alignedat}{3}
&\llangle\boldsymbol{P}_{h},\mathbf{grad}\,\boldsymbol{\Upsilon}_{h}\rrangle=\llangle \boldsymbol{B},\boldsymbol{\Upsilon}_{h} \rrangle + \int_{\Gamma_{2}}\overline{\boldsymbol{T}}\boldsymbol{\cdot}\boldsymbol{\Upsilon}_{h} dA, &\quad &\forall \boldsymbol{\Upsilon}_{h}\in [V^{1}_{h,1}]^{n},\\
&\llangle\mathbf{grad}\,\boldsymbol{U}_{h},\boldsymbol{\lambda}_{h}\rrangle - \llangle \boldsymbol{K}_{h},\boldsymbol{\lambda}_{h}\rrangle = 0,& &\forall \boldsymbol{\lambda}_{h}\in V^{\mathbf{c}}_{h},\\
&\llangle \mathbb{P}(\boldsymbol{K}_{h}),\boldsymbol{\pi}_{h}\rrangle - \llangle \boldsymbol{P}_{h}, \boldsymbol{\pi}_{h} \rrangle - \llangle p_{h}\,(\boldsymbol{I} +\boldsymbol{K}_{h})^{-T}, \boldsymbol{\pi}_{h} \rrangle  = 0, & &\forall \boldsymbol{\pi}_{h}\in V^{\mathbf{d}}_{h}, \\
&\llangle \det(\boldsymbol{I} + \boldsymbol{K}_{h}), q_{h} \rrangle = \llangle 1,q_{h}\rrangle, & &\forall q_{h}\in V^{\mathrm{D}}_{h}.  
\end{alignedat}
\end{equation}
\end{minipage}
\bigskip

We employ simplicial finite elements shown in Figure \ref{FEO012} to generate the above finite element spaces in $2$D and $3$D. More specifically, we use the standard simplicial Lagrange elements and discontinuous elements respectively for $[V^{1}_{h}]^{n}$ and $V^{\mathrm{D}}_{h}$. As mentioned earlier, the definitions of $\mathbf{curl}$ and $\mathbf{div}$ for second-order tensors imply that in $\mathbb{R}^{n}$, the tensorial spaces $H^{\mathbf{c}}(\mathcal{B})$ and $H^{\mathbf{d}}(\mathcal{B})$ are equivalent to $n$-copies of the associated vectorial spaces in the sense that each row of second-order tensors of the $H^{\mathbf{c}}$- and $H^{\mathbf{d}}$-classes can be considered as a vector field beloning to the vectorial spaces $H^{c}(\mathcal{B})$ and $H^{d}(\mathcal{B})$ corresponding to the standard $\mathrm{curl}$ and $\mathrm{div}$ operators of vector fields, respectively. We use the N\'{e}d\'{e}lec elements of the first and the second kinds to generate $V^{\mathbf{c}}_{h}$ and the Raviart-Thomas elements and the Brezzi-Douglas-Marini elements to generate $V^{\mathbf{d}}_{h}$.

\section{Newtons' Iterations}
We employ Newton's method to solve \eqref{IncompFEMs}. Here we mention this approach in some details since, as will be discussed in the following section, it is also useful for stability analysis. Let $\bar{Z}=[H^{1}(\mathcal{B})]^{n}\times H^{\mathbf{c}}(\mathcal{B}) \times H^{\mathbf{d}}(\mathcal{B})\times L^{2}(\mathcal{B})$ and $Z=[H^{1}_{1}(\mathcal{B})]^{n}\times H^{\mathbf{c}}(\mathcal{B}) \times H^{\mathbf{d}}(\mathcal{B})\times L^{2}(\mathcal{B})$. Also let $u\in \bar{Z}$ and $y,z\in Z$, with $u=(\boldsymbol{U},\boldsymbol{K},\boldsymbol{P},p)$, $y=(\boldsymbol{\Upsilon},\boldsymbol{\lambda},\boldsymbol{\pi},q)$, and $z=(\boldsymbol{V},\boldsymbol{M},\boldsymbol{Q},r)$. The nonlinear problem \eqref{IncompWeakForm} can be written as: Find $u\in \bar{Z}$ such that
\begin{equation}\label{Abs_NLProblem}
\begin{alignedat}{3}
\langle H(u), y \rangle &= \llangle\boldsymbol{P},\mathbf{grad}\,\boldsymbol{\Upsilon}\rrangle + \llangle\mathbf{grad}\,\boldsymbol{U},\boldsymbol{\lambda}\rrangle - \llangle \boldsymbol{K},\boldsymbol{\lambda}\rrangle\\
 &+ \llangle \mathbb{P}(\boldsymbol{K}),\boldsymbol{\pi}\rrangle - \llangle \boldsymbol{P}, \boldsymbol{\pi} \rrangle - \llangle p\,(\boldsymbol{I} +\boldsymbol{K})^{-T}, \boldsymbol{\pi} \rrangle \\
&+ \llangle \det(\boldsymbol{I} + \boldsymbol{K}), q \rrangle - \llangle \boldsymbol{B},\boldsymbol{\Upsilon} \rrangle - \int_{\Gamma_{2}}\overline{\boldsymbol{T}}\boldsymbol{\cdot}\boldsymbol{\Upsilon} dA - \llangle 1,q\rrangle = 0, \quad \forall y \in Z. 
\end{alignedat}
\end{equation}
The $i$-th Newton's iteration for solving the above nonlinear equation reads: Given $u^{i}\in \bar{Z}$, find $z\in Z$ by solving a linear problem and let $u^{i+1}=u^{i}+z$. This linear problem is obtained by linearizing the above equation and can be stated as: Find $z\in Z$ such that 
\begin{align}\label{ith_NewItr}
b(z,y;u^{i}) = -\langle H(u^{i}),y \rangle, \quad \forall y \in Z,
\end{align}  
where the bilinear form $b(\cdot,\cdot;u^{i})$ is given by
\begin{equation}\label{BiLin_Newt}
\begin{aligned}
b(z, y;u^{i}) &= \llangle\boldsymbol{Q},\mathbf{grad}\,\boldsymbol{\Upsilon}\rrangle + \llangle\mathbf{grad}\,\boldsymbol{V},\boldsymbol{\lambda}\rrangle - \llangle \boldsymbol{M},\boldsymbol{\lambda}\rrangle + \llangle \mathsf{A}(\boldsymbol{K}^{i})\boldsymbol{:}\boldsymbol{M},\boldsymbol{\pi}\rrangle - \llangle \boldsymbol{Q}, \boldsymbol{\pi} \rrangle \\
&+ \llangle p^{i} (\boldsymbol{I} + \boldsymbol{K}^{i})^{-T} \boldsymbol{M}^{T}(\boldsymbol{I} + \boldsymbol{K}^{i})^{-T}, \boldsymbol{\pi} \rrangle - \llangle r\,(\boldsymbol{I} + \boldsymbol{K}^{i})^{-T}, \boldsymbol{\pi} \rrangle \\
&+  \llangle \det(\boldsymbol{I} + \boldsymbol{K}^{i}) \,\, \mathrm{tr}\!\left[ (\boldsymbol{I} + \boldsymbol{K}^{i})^{-1} \boldsymbol{M} \right] , q\rrangle, 
\end{aligned}
\end{equation}  
where $\mathsf{A}(\boldsymbol{K})$ is the elasticity tensor in terms of the displacement gradient and $(\mathsf{A}(\boldsymbol{K})\boldsymbol{:}\boldsymbol{M})^{IJ}:=A^{IJRS}M^{RS}$. By discretizing the linear problem \eqref{ith_NewItr} using the finite element spaces of the previous section, one obtains the following Newton's iteration for the finite element method \eqref{IncompFEMs}:  

\bigskip
\noindent{\textit{Let $\bar{Z}_{h} = [V^{1}_{h}]^{n}\times V^{\mathbf{c}}_{h} \times V^{\mathbf{d}}_{h}\times V^{\mathrm{D}}_{h}$, and $Z_{h} = [V^{1}_{h,1}]^{n}\times V^{\mathbf{c}}_{h} \times V^{\mathbf{d}}_{h}\times V^{\mathrm{D}}_{h}$. Given $u^{i}_{h}=(\boldsymbol{U}^{i}_{h},\boldsymbol{K}^{i}_{h},\boldsymbol{P}^{i}_{h},p^{i}_{h})\in \bar{Z}_{h} $, find $z_{h}=(\boldsymbol{V}_{h},\boldsymbol{M}_{h},\boldsymbol{Q}_{h},r_{h})\in Z_{h} $ and let $u^{i+1}_{h} = u^{i}_{h} + z_{h}$, where $z_{h}$ is obtained by solving the linear finite element method}}
\begin{equation}\label{NewtIt_IncompFEMs}
\begin{alignedat}{3}
&\llangle\boldsymbol{Q}_{h},\mathbf{grad}\,\boldsymbol{\Upsilon}_{h}\rrangle = -\llangle\boldsymbol{P}^{i}_{h},\mathbf{grad}\,\boldsymbol{\Upsilon}_{h}\rrangle+\llangle \boldsymbol{B},\boldsymbol{\Upsilon}_{h} \rrangle + \int_{\Gamma_{2}}\overline{\boldsymbol{T}}\boldsymbol{\cdot}\boldsymbol{\Upsilon}_{h} dA, &\quad &\forall \boldsymbol{\Upsilon}_{h}\in [V^{1}_{h,1}]^{n},\\
&\llangle\mathbf{grad}\,\boldsymbol{V}_{h},\boldsymbol{\lambda}_{h}\rrangle - \llangle \boldsymbol{M}_{h},\boldsymbol{\lambda}_{h}\rrangle = -\llangle\mathbf{grad}\,\boldsymbol{U}^{i}_{h},\boldsymbol{\lambda}_{h}\rrangle + \llangle \boldsymbol{K}^{i}_{h},\boldsymbol{\lambda}_{h}\rrangle, & & \forall \boldsymbol{\lambda}_{h}\in V^{\mathbf{c}}_{h},\\
&\llangle \mathsf{A}(\boldsymbol{K}^{i}_{h})\boldsymbol{:}\boldsymbol{M}_{h},\boldsymbol{\pi}_{h}\rrangle - \llangle \boldsymbol{Q}_{h}, \boldsymbol{\pi}_{h} \rrangle + \llangle p^{i}_{h} (\boldsymbol{I} + \boldsymbol{K}^{i}_{h})^{-T} \boldsymbol{M}^{T}_{h}(\boldsymbol{I} + \boldsymbol{K}^{i}_{h})^{-T}, \boldsymbol{\pi}_{h} \rrangle &&\\
& \quad \quad - \llangle r_{h}(\boldsymbol{I} + \boldsymbol{K}^{i}_{h})^{-T}, \boldsymbol{\pi}_{h} \rrangle =  -\llangle \mathbb{P}(\boldsymbol{K}^{i}_{h}),\boldsymbol{\pi}_{h}\rrangle + \llangle \boldsymbol{P}^{i}_{h}, \boldsymbol{\pi}_{h} \rrangle & &\\
& \hspace{2in} + \llangle p^{i}_{h}\,(\boldsymbol{I} +\boldsymbol{K}^{i}_{h})^{-T}, \boldsymbol{\pi}_{h} \rrangle, & & \forall \boldsymbol{\pi}_{h}\in V^{\mathbf{d}}_{h}, \\
&\llangle \det(\boldsymbol{I} + \boldsymbol{K}^{i}_{h}) \,\, \mathrm{tr}\!\left[ (\boldsymbol{I} + \boldsymbol{K}^{i}_{h})^{-1} \boldsymbol{M}_{h} \right] , q_{h}\rrangle =
-\llangle \det(\boldsymbol{I} + \boldsymbol{K}^{i}_{h}), q_{h} \rrangle + \llangle 1, q_{h} \rrangle, && \forall q_{h}\in V^{\mathrm{D}}_{h}.  
\end{alignedat}
\end{equation}
\bigskip

Suppose 
\begin{align}\label{Dim_FESp}
n_{1} = \dim [V^{1}_{h,1}]^{n},\quad n_{\mathbf{c}} = \dim V^{\mathbf{c}}_{h}, \quad n_{\mathbf{d}} = \dim V^{\mathbf{d}}_{h}, \quad n_{\mathrm{D}} = \dim V^{\mathrm{D}}_{h},
\end{align}
and let $n_{t}=n_{1}+n_{\mathbf{c}}+n_{\mathbf{d}}+n_{\mathrm{D}}$, denote the total number of degrees of freedom. Then, it is straightforward to see that the stiffness matrix of the above linear problem is of the form  
\begin{align}\label{LinStiffM}
\mathbb{S}_{n_{t}\times n_{t}}=\left[ \arraycolsep=1.1pt\def\arraystretch{1.2} \begin{array}{c;{2pt/2pt}c;{2pt/2pt}c;{2pt/2pt}c} \mathbf{0} & \mathbf{0} & \mathbb{S}^{1\mathbf{d}}_{n_{1}\times n_{\mathbf{d}}} & \mathbf{0} \\ \hdashline[2pt/2pt] 
\mathbb{S}^{\mathbf{c}1}_{n_{\mathbf{c}}\times n_{1}} & \mathbb{S}^{\mathbf{c}\mathbf{c}}_{n_{\mathbf{c}}\times n_{\mathbf{c}}} &\mathbf{0} & \mathbf{0} \\ \hdashline[2pt/2pt]
\mathbf{0} & \mathbb{S}^{\mathbf{d}\mathbf{c}}_{n_{\mathbf{d}}\times n_{\mathbf{c}}} & \mathbb{S}^{\mathbf{d}\mathbf{d}}_{n_{\mathbf{d}}\times n_{\mathbf{d}}} & \mathbb{S}^{\mathbf{d}\mathrm{D}}_{n_{\mathbf{d}}\times n_{\mathrm{D}}} \\ \hdashline[2pt/2pt]
\mathbf{0} & \mathbb{S}^{\mathrm{D}\mathbf{c}}_{n_{\mathrm{D}}\times n_{\mathbf{c}}} & \mathbf{0} & \mathbf{0} \end{array}   \right].		
\end{align}


\section{Stability Analysis}\label{Sec_StAnalysis}
The stability of the nonlinear finite element method \eqref{IncompFEMs} can be studied by using the general theory for the Galerkin approximation of nonlinear problems discussed in \citep{PoRa1994,CaRa1997}. Let $u=(\boldsymbol{U},\boldsymbol{K},\boldsymbol{P},p)$ be a regular solution of \eqref{IncompWeakForm}, or equivalently \eqref{Abs_NLProblem}, in the sense that the derivative of the nonlinear mapping $H$ defined in \eqref{Abs_NLProblem} is nonsingular at $u$. This derivative can be expressed in terms of the bilinear form $b(\cdot,\cdot;u)$ introduced in \eqref{BiLin_Newt}. The theory of \citep{PoRa1994,CaRa1997} states that in a neighborhood of $u$ and for sufficiently refined meshes with the maximum element diameter $h>0$, under some additional mild conditions, the nonlinear finite element method \eqref{IncompFEMs} has a unique solution $u_{h}=(\boldsymbol{U}_{h},\boldsymbol{K}_{h},\boldsymbol{P}_{h},p_{h})$ that converges to $u$ as $h\rightarrow 0$ if there exists a mesh-independent number $\alpha>0$ such that
\begin{equation}\label{infsup_suf}
\underset{y_{h}\in Z_{h}}{\inf}\, \underset{z_{h}\in Z_{h}}{\sup} \frac{b(z_{h},y_{h};u)}{\|z_{h}\|_{Z} \|y_{h}\|_{Z} } \geq \alpha>0,
\end{equation}
where $Z_{h}$ is the mixed finite element space for the Newton's iteration \eqref{NewtIt_IncompFEMs} and $\|\cdot\|_{Z}$ is the norm of $Z$. Thus, the inf-sup condition \eqref{infsup_suf} is a sufficient stability condition.

The condition \eqref{infsup_suf} can be interpreted as a stability condition for the linear finite element method of the Newton's iteration \eqref{NewtIt_IncompFEMs} as well. More specifically, suppose that at the $(i-1)$-th Newton's iteration, we have $u \approx u^{i}_{h}$. Then, the condition \eqref{infsup_suf} implies that the linear problem \eqref{NewtIt_IncompFEMs} at the the $i$-th iteration has a unique solution and the stiffness matrix $\mathbb{S}$ introduced in \eqref{LinStiffM} is non-singular. One can also approximate the lower bound $\alpha$ of \eqref{infsup_suf} by using a matrix associated to $\mathbb{S}$ \citep[Proposition 3.4.5]{BoBrFo2013}. In particular, $\alpha$ can be approximated by the smallest singular value of the matrix $\mathbb{M}\mathbb{S}\mathbb{M}$, where $\mathbb{M}$ is the symmetric, positive definite matrix associated to the norm $\|\cdot\|_{Z}$.         

The non-singularity of the stiffness matrix $\mathbb{S}$ implies that the submatrices 
\begin{alignat}{3}\label{SubMat}
&\mathbb{S}^{1\mathbf{d}}_{n_{1}\times n_{\mathbf{d}}}, \quad \mathbb{S}^{\mathrm{D}\mathbf{c}}_{n_{\mathrm{D}}\times n_{\mathbf{c}}},  
& \quad &\mathbb{B}=\left[ \arraycolsep=1.1pt\def\arraystretch{1.2} \begin{array}{c;{2pt/2pt}c;{2pt/2pt}c} \mathbf{0} & \mathbf{0} & \mathbb{S}^{1\mathbf{d}}_{n_{1}\times n_{\mathbf{d}}} \\ \hdashline[2pt/2pt] 
\mathbb{S}^{\mathbf{c}1}_{n_{\mathbf{c}}\times n_{1}} & \mathbb{S}^{\mathbf{c}\mathbf{c}}_{n_{\mathbf{c}}\times n_{\mathbf{c}}} &\mathbf{0} \end{array}   \right],  
& \quad & \mathbb{C}=\left[ \arraycolsep=1.1pt\def\arraystretch{1.2} \begin{array}{c;{2pt/2pt}c;{2pt/2pt}c} \mathbf{0} & \mathbb{S}^{1\mathbf{d}}_{n_{1}\times n_{\mathbf{d}}} & \mathbf{0} \\ \hdashline[2pt/2pt] 
\mathbb{S}^{\mathbf{d}\mathbf{c}}_{n_{\mathbf{d}}\times n_{\mathbf{c}}} & \mathbb{S}^{\mathbf{d}\mathbf{d}}_{n_{\mathbf{d}}\times n_{\mathbf{d}}} & \mathbb{S}^{\mathbf{d}\mathrm{D}}_{n_{\mathbf{d}}\times n_{\mathrm{D}}} \end{array}   \right],\\
& \mathbb{D}=\left[ \arraycolsep=1.1pt\def\arraystretch{1.2} \begin{array}{c;{2pt/2pt}c} \mathbb{S}^{\mathbf{c}1}_{n_{\mathbf{c}}\times n_{1}} & \mathbb{S}^{\mathbf{c}\mathbf{c}}_{n_{\mathbf{c}}\times n_{\mathbf{c}}} \\ \hdashline[2pt/2pt]
\mathbf{0} & \mathbb{S}^{\mathrm{D}\mathbf{c}}_{n_{\mathrm{D}}\times n_{\mathbf{c}}}\end{array}   \right], 
& &\mathbb{E}=\left[ \arraycolsep=1.1pt\def\arraystretch{1.2} \begin{array}{c;{2pt/2pt}c;{2pt/2pt}c} \mathbb{S}^{\mathbf{d}\mathbf{c}}_{n_{\mathbf{d}}\times n_{\mathbf{c}}} & \mathbb{S}^{\mathbf{d}\mathbf{d}}_{n_{\mathbf{d}}\times n_{\mathbf{d}}} & \mathbb{S}^{\mathbf{d}\mathrm{D}}_{n_{\mathbf{d}}\times n_{\mathrm{D}}} \\ \hdashline[2pt/2pt]
\mathbb{S}^{\mathrm{D}\mathbf{c}}_{n_{\mathrm{D}}\times n_{\mathbf{c}}} & \mathbf{0} & \mathbf{0} \end{array} \right], &&		
\end{alignat}
must be full rank. This result can be stated by using some inf-sup conditions associated to suitable bilinear forms induced by \eqref{NewtIt_IncompFEMs}. The upshot can be stated as follows:
\begin{itemize}
\item $\mathbb{S}^{1\mathbf{d}}$ is full rank if and only if there exists $\alpha_{h}>0$ such that
\begin{align}\label{infsup_S1d}
\underset{\boldsymbol{\Upsilon}_{h}\in [V^{1}_{h,1}]^{n}}{\inf}\, \underset{\boldsymbol{Q}_{h}\in V^{\mathbf{d}}_{h}}{\sup} \frac{\llangle\boldsymbol{Q}_{h},\mathbf{grad}\,\boldsymbol{\Upsilon}_{h}\rrangle}{\|\boldsymbol{Q}_{h}\|_{\mathbf{d}}\, \|\boldsymbol{\Upsilon}_{h}\|_{1} } \geq \alpha_{h};
\end{align}
\item $\mathbb{S}^{\mathrm{D}\mathbf{c}}$ is full rank if and only if there exists $\alpha_{h}>0$ such that
\begin{align}\label{infsup_SDc}
\underset{q_{h}\in V^{\mathrm{D}}_{h}}{\inf}\, \underset{\boldsymbol{M}_{h}\in V^{\mathbf{c}}_{h}}{\sup} \frac{\llangle \det(\boldsymbol{I} + \boldsymbol{K}^{i}_{h}) \,\, \mathrm{tr}\!\left[ (\boldsymbol{I} + \boldsymbol{K}^{i}_{h})^{-1} \boldsymbol{M}_{h} \right] , q_{h}\rrangle}{\|\boldsymbol{M}_{h}\|\, \|q_{h}\| } \geq \alpha_{h};
\end{align}
\item $\mathbb{B}$ is full rank if and only if there exists $\alpha_{h}>0$ such that
\begin{align}\label{infsup_B}
\underset{(\boldsymbol{\Upsilon}_{h},\boldsymbol{\lambda}_{h})\in [V^{1}_{h,1}]^{n} \times V^{\mathbf{c}}_{h}}{\inf}\, \underset{(\boldsymbol{V}_{h},\boldsymbol{M}_{h},\boldsymbol{Q}_{h})\in [V^{1}_{h,1}]^{n}\times V^{\mathbf{c}}_{h} \times V^{\mathbf{d}}_{h}}{\sup} \frac{b^{\mathbb{B}}\big((\boldsymbol{V}_{h},\boldsymbol{M}_{h},\boldsymbol{Q}_{h}),(\boldsymbol{\Upsilon}_{h},\boldsymbol{\lambda}_{h})\big)}{\|(\boldsymbol{V}_{h},\boldsymbol{M}_{h},\boldsymbol{Q}_{h})\|\, \|(\boldsymbol{\Upsilon}_{h},\boldsymbol{\lambda}_{h})\| } \geq \alpha_{h},
\end{align}
where
\begin{align}
b^{\mathbb{B}}\big((\boldsymbol{V}_{h},\boldsymbol{M}_{h},\boldsymbol{Q}_{h}),(\boldsymbol{\Upsilon}_{h},\boldsymbol{\lambda}_{h})\big) = \llangle\boldsymbol{Q}_{h},\mathbf{grad}\,\boldsymbol{\Upsilon}_{h}\rrangle + \llangle\mathbf{grad}\,\boldsymbol{V}_{h},\boldsymbol{\lambda}_{h}\rrangle - \llangle \boldsymbol{M}_{h},\boldsymbol{\lambda}_{h}\rrangle;
\end{align}
\item $\mathbb{C}$ is full rank if and only if there exists $\alpha_{h}>0$ such that
\begin{align}\label{infsup_C}
\underset{(\boldsymbol{\Upsilon}_{h},\boldsymbol{\pi}_{h})\in [V^{1}_{h,1}]^{n} \times V^{\mathbf{d}}_{h}}{\inf}\, \underset{(\boldsymbol{M}_{h},\boldsymbol{Q}_{h},r_{h})\in V^{\mathbf{c}}_{h}\times V^{\mathbf{d}}_{h} \times V^{\mathrm{D}}_{h} }{\sup} \frac{b^{\mathbb{C}}\big((\boldsymbol{M}_{h},\boldsymbol{Q}_{h},r_{h}),(\boldsymbol{\Upsilon}_{h},\boldsymbol{\pi}_{h}); \boldsymbol{K}^{i}_{h}, p^{i}_{h}\big)}{\|(\boldsymbol{M}_{h},\boldsymbol{Q}_{h},r_{h})\|\, \|(\boldsymbol{\Upsilon}_{h},\boldsymbol{\pi}_{h})\| } \geq \alpha_{h},
\end{align}
where
\begin{alignat}{3}
b^{\mathbb{C}}\big((\boldsymbol{M}_{h},&\boldsymbol{Q}_{h},r_{h}),(\boldsymbol{\Upsilon}_{h},\boldsymbol{\pi}_{h}); \boldsymbol{K}^{i}_{h}, p^{i}_{h}\big) = \llangle\boldsymbol{Q}_{h},\mathbf{grad}\,\boldsymbol{\Upsilon}_{h}\rrangle + \llangle \mathsf{A}(\boldsymbol{K}^{i}_{h})\boldsymbol{:}\boldsymbol{M}_{h},\boldsymbol{\pi}_{h}\rrangle \\
& - \llangle \boldsymbol{Q}_{h}, \boldsymbol{\pi}_{h} \rrangle + \llangle p^{i}_{h} (\boldsymbol{I} + \boldsymbol{K}^{i}_{h})^{-T} \boldsymbol{M}^{T}_{h}(\boldsymbol{I} + \boldsymbol{K}^{i}_{h})^{-T}, \boldsymbol{\pi}_{h} \rrangle - \llangle r_{h}(\boldsymbol{I} + \boldsymbol{K}^{i}_{h})^{-T}, \boldsymbol{\pi}_{h} \rrangle;
\end{alignat}
\item $\mathbb{D}$ is full rank if and only if there exists $\alpha_{h}>0$ such that
\begin{align}\label{infsup_D}
\underset{(\boldsymbol{\lambda}_{h},q_{h})\in V^{\mathbf{c}}_{h}\times V^{\mathrm{D}}_{h}}{\inf}\, \underset{(\boldsymbol{V}_{h},\boldsymbol{M}_{h})\in [V^{1}_{h,1}]^{n} \times V^{\mathbf{c}}_{h}}{\sup} \frac{b^{\mathbb{D}}\big((\boldsymbol{V}_{h},\boldsymbol{M}_{h}),(\boldsymbol{\lambda}_{h},q_{h}); \boldsymbol{K}^{i}_{h}\big)}{\|(\boldsymbol{V}_{h},\boldsymbol{M}_{h})\|\, \|(\boldsymbol{\lambda}_{h},q_{h})\| } \geq \alpha_{h},
\end{align}
where
\begin{alignat}{3}
b^{\mathbb{D}}\big((\boldsymbol{V}_{h},\boldsymbol{M}_{h}),(\boldsymbol{\lambda}_{h},q_{h}); \boldsymbol{K}^{i}_{h}\big) &= \llangle\mathbf{grad}\,\boldsymbol{V}_{h},\boldsymbol{\lambda}_{h}\rrangle - \llangle \boldsymbol{M}_{h},\boldsymbol{\lambda}_{h}\rrangle \\
 &+ \llangle \det(\boldsymbol{I} + \boldsymbol{K}^{i}_{h}) \,\, \mathrm{tr}\!\left[ (\boldsymbol{I} + \boldsymbol{K}^{i}_{h})^{-1} \boldsymbol{M}_{h} \right] , q_{h}\rrangle;
\end{alignat}
\item $\mathbb{E}$ is full rank if and only if there exists $\alpha_{h}>0$ such that
\begin{align}\label{infsup_E}
\underset{(\boldsymbol{\pi}_{h},q_{h})\in V^{\mathbf{d}}_{h} \times V^{\mathrm{D}}_{h}}{\inf}\, \underset{(\boldsymbol{M}_{h},\boldsymbol{Q}_{h},r_{h})\in V^{\mathbf{c}}_{h}\times V^{\mathbf{d}}_{h} \times V^{\mathrm{D}}_{h} }{\sup} \frac{b^{\mathbb{E}}\big((\boldsymbol{M}_{h},\boldsymbol{Q}_{h},r_{h}),(\boldsymbol{\pi}_{h},q_{h}); \boldsymbol{K}^{i}_{h}, p^{i}_{h}\big)}{\|(\boldsymbol{M}_{h},\boldsymbol{Q}_{h},r_{h})\|\, \|(\boldsymbol{\pi}_{h},q_{h})\| } \geq \alpha_{h},
\end{align}
where
\begin{alignat}{3}
b^{\mathbb{E}}\big((\boldsymbol{M}_{h},\boldsymbol{Q}_{h},r_{h}),&(\boldsymbol{\pi}_{h},q_{h}); \boldsymbol{K}^{i}_{h}, p^{i}_{h}\big) = \llangle \mathsf{A}(\boldsymbol{K}^{i}_{h})\boldsymbol{:}\boldsymbol{M}_{h},\boldsymbol{\pi}_{h}\rrangle - \llangle \boldsymbol{Q}_{h}, \boldsymbol{\pi}_{h} \rrangle \\
& + \llangle p^{i}_{h} (\boldsymbol{I} + \boldsymbol{K}^{i}_{h})^{-T} \boldsymbol{M}^{T}_{h}(\boldsymbol{I} + \boldsymbol{K}^{i}_{h})^{-T}, \boldsymbol{\pi}_{h} \rrangle - \llangle r_{h}(\boldsymbol{I} + \boldsymbol{K}^{i}_{h})^{-T}, \boldsymbol{\pi}_{h} \rrangle \\
& + \llangle \det(\boldsymbol{I} + \boldsymbol{K}^{i}_{h}) \,\, \mathrm{tr}\!\left[ (\boldsymbol{I} + \boldsymbol{K}^{i}_{h})^{-1} \boldsymbol{M}_{h} \right] , q_{h}\rrangle.
\end{alignat}
\end{itemize}

A simple inspection of the number of rows and columns of the above matrices together with the rank-nullity theorem yield the following necessary condition for the validity of \eqref{infsup_suf} and the well-posedness of \eqref{NewtIt_IncompFEMs}, where we use the notation of \eqref{Dim_FESp}.

\begin{thm}\label{NecInequal} The inf-sup condition \eqref{infsup_suf} does not hold and the Newton's iteration \eqref{NewtIt_IncompFEMs} does not admit a unique solution if at least one of the following inequalities holds: (i) $n_{1}> n_{\mathbf{d}}$; (ii) $n_{\mathrm{D}} > n_{\mathbf{c}}$; (iii) $n_{1} > n_{\mathbf{c}} + n_{\mathrm{D}}$; (iv) $n_{\mathrm{D}} > n_{1}$.
\end{thm}

By using this result, one can specify some unstable combinations of finite elements. For example, consider the $2$D  elements of Figure \ref{FEO012} and suppose $N_{v}$, $N_{ed}$, and $N_{el}$ are respectively the number of vertices, edges, and elements of a $2$D simplicial mesh. One can write the relations
\begin{align}\label{2DMesh_Rela}
N_{el} - N_{ed} + N_{v} = 1 - I, \text{ and } 2N_{ed}-N^{\partial}_{ed} = 3 N_{el},
\end{align}
where $I$ is the number of holes and $N^{\partial}_{ed}$ is the number boundary edges \citep{ErnGuermond2004}. By using the notation of Figure \ref{FEO012}, we obtain the following relations on a $2$D simplicial mesh:
\begin{alignat}{3}\label{FESp_dim}
&n_{\mathrm{D}} = \left\{ \begin{array}{ll} 3N_{el}, &\text{ $\mathrm{D1}$ element,} \\
6N_{el}, &\text{ $\mathrm{D2}$ element},\end{array} \right.
\quad & &n_{1} = \left\{ \begin{array}{ll} 2N_{v}, &\text{ $\mathrm{L1}$ element,} \\
2(N_{v}+N_{ed}), &\text{ $\mathrm{L2}$ element},\end{array} \right.\\
&n_{\mathbf{c}} = 2N_{ed}, \text{ $\mathrm{N11}$ element}, & &n_{\mathbf{d}} = 2N_{ed}, \text{ $\mathrm{R1}$ element}.
\end{alignat} 
A simple varification of the inequalities of Theorem \ref{NecInequal} by using the relations \eqref{2DMesh_Rela} and \eqref{FESp_dim} yields the following result.

\begin{cor}\label{2D_instFEs} Let $\mathrm{FE}_{1}$, $\mathrm{FE}_{\mathbf{c}}$, $\mathrm{FE}_{\mathbf{d}}$, and $\mathrm{FE}_{\mathrm{D}}$ be respectively arbitrary finite element spaces for $H^{1}$, $H^{\mathbf{c}}$, $H^{\mathbf{d}}$, and $L^{2}$ spaces. In $2$D, the choices $(\mathrm{L2}, \mathrm{FE}_{\mathbf{c}}, \mathrm{R1}, \mathrm{FE}_{\mathrm{D}})$, $(\mathrm{FE}_{1},\mathrm{N11}, \mathrm{FE}_{\mathbf{d}}, \mathrm{D2})$, $(\mathrm{L1}, \mathrm{FE}_{\mathbf{c}}, \mathrm{FE}_{\mathbf{d}}, \mathrm{D1})$, and $(\mathrm{L1}, \mathrm{FE}_{\mathbf{c}}, \mathrm{FE}_{\mathbf{d}}, \mathrm{D2})$ for the finite element methods \eqref{NewtIt_IncompFEMs} violate at least one of the inequalities of Theorem \ref{NecInequal} and therefore, lead to unstable finite element methods.
\end{cor}

\begin{figure}[t]
\begin{center}
\includegraphics[scale=.55,angle=0]{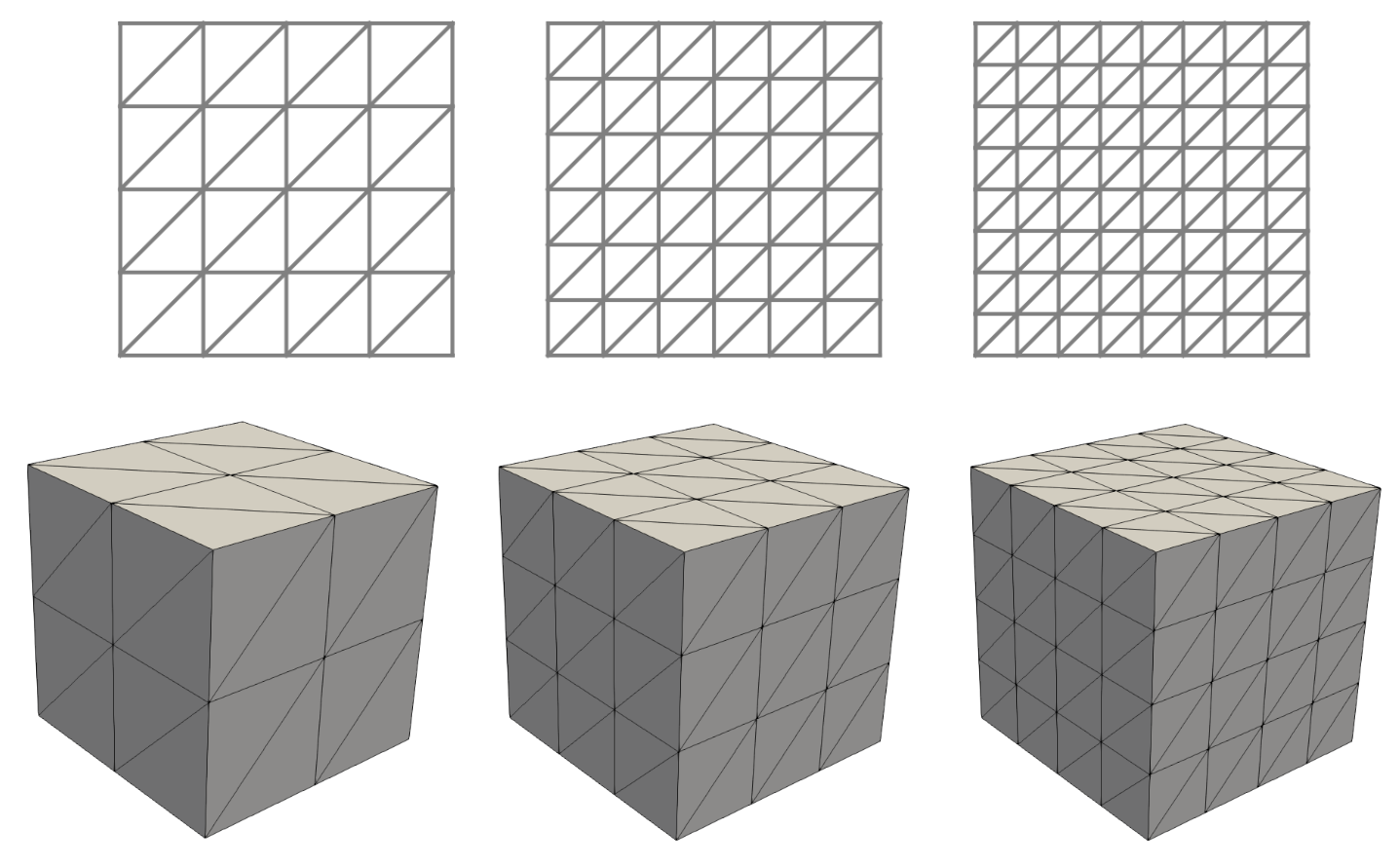}
\end{center}
\vspace*{-0.3in}
\caption{\footnotesize Simplicial meshes of a unit square and a unit cube used in the numerical examples.} 
\label{Meshes}
\end{figure}

\section{Numerical Results}
Some numerical examples are discussed in this section to study the stability and the convergence of the mixed finite element methods \eqref{IncompFEMs}. We consider the simplicial finite elements of Figure \ref{FEO012} and use FEniCS \citep{loMaWe2012} to implement these numerical examples. The underlying domains are assumed to be a unit square in $2$D and a unit cube in $3$D with simplicial meshes shown in Figure \ref{Meshes}. We consider incompressible Neo-Hookean materials with the stored energy function $W(\boldsymbol{F})=\frac{\mu}{2}(\mathrm{tr}\boldsymbol{F^{T}}\boldsymbol{F}-3)$,  $\mu>0$, and the constitutive equation
\begin{equation}\label{IncompNeoHook} 
\boldsymbol{P}(\boldsymbol{F})= \mu\boldsymbol{F} - p\boldsymbol{F}^{-T}.
\end{equation}
The term $\mathsf{A}(\boldsymbol{K})\boldsymbol{:}\boldsymbol{M}$ of the bilinear form \eqref{BiLin_Newt} then simply reads $
\mathsf{A}(\boldsymbol{K})\boldsymbol{:}\boldsymbol{M}= \mu \boldsymbol{M}$.

\paragraph{Notation.}To refer to any choice of the elements of Figure \ref{FEO012} for the discretization of \eqref{IncompWeakForm}, we use the abbreviated names shown in Figure \ref{FEO012}. For example, $\mathrm{L1N12R2D0}$ indicates the choice of the Lagrange element of degree $1$, the degree $2$ N\'{e}d\'{e}lec element of the first kind, the Raviart-Thomas element of degree $2$, and the discontinuous element of degree $0$ respectively for the discretization of $(\boldsymbol{U},\boldsymbol{K},\boldsymbol{P},p)$.

\subsection{Stability Study}
\begin{figure}[t]
\begin{center}
\includegraphics[scale=.52,angle=0]{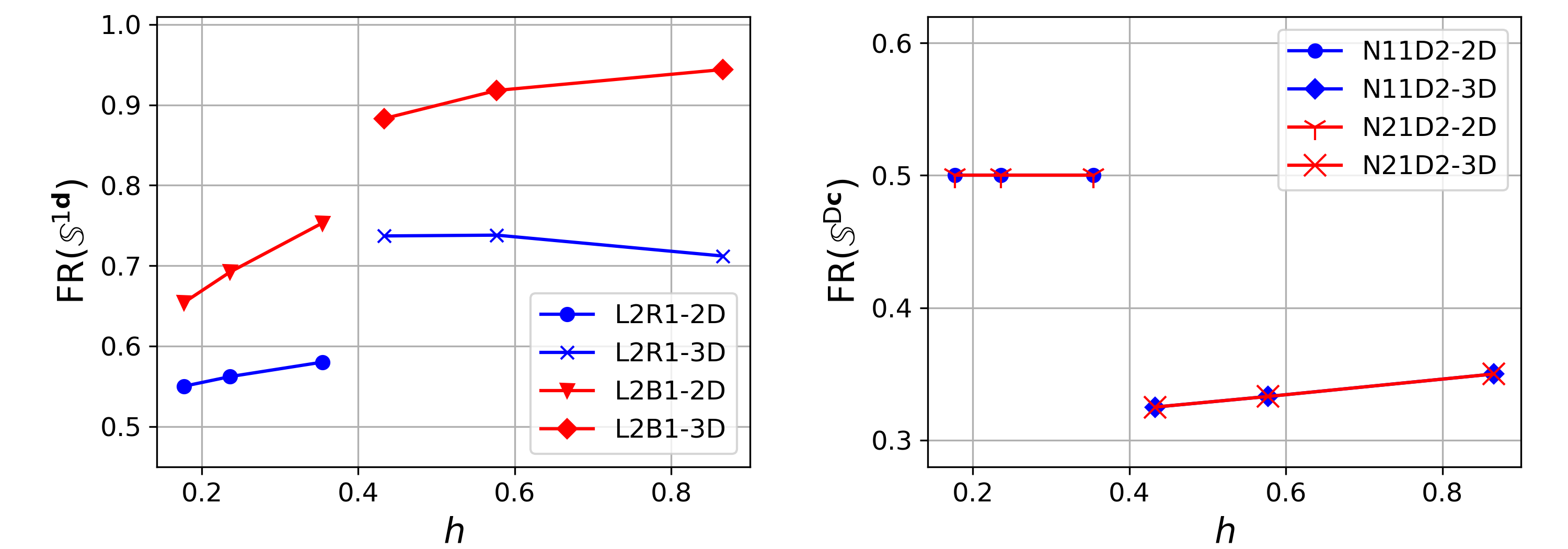}
\includegraphics[scale=.52,angle=0]{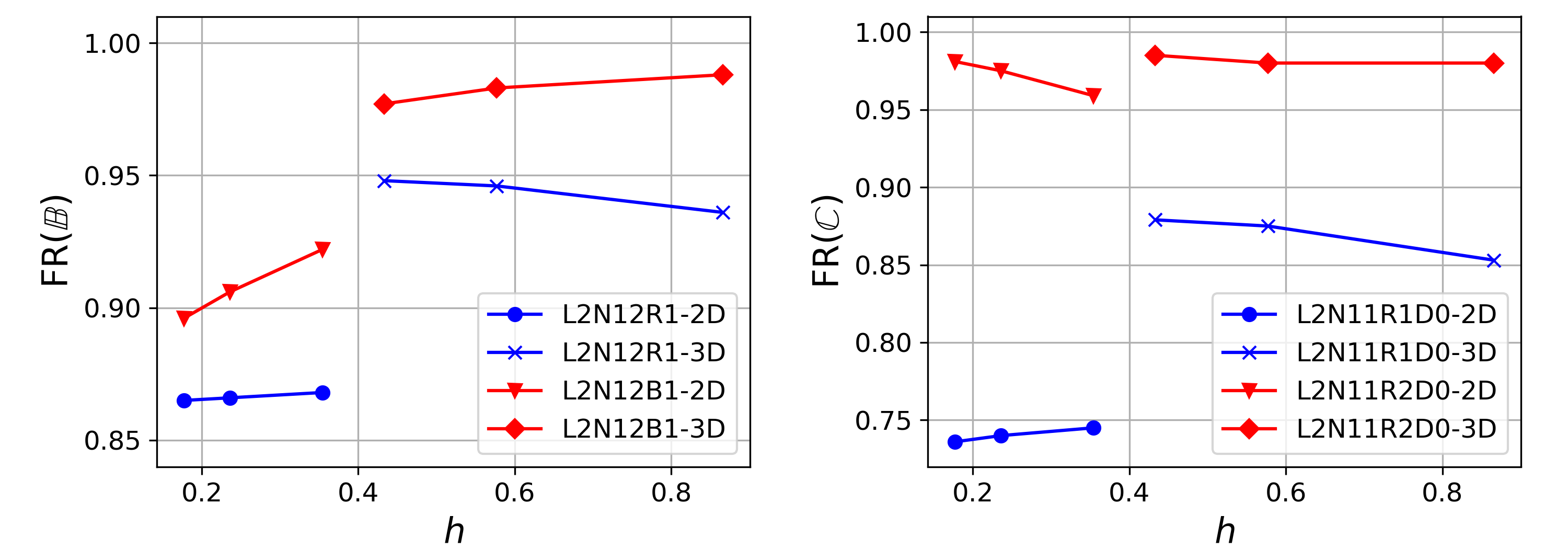}
\includegraphics[scale=.52,angle=0]{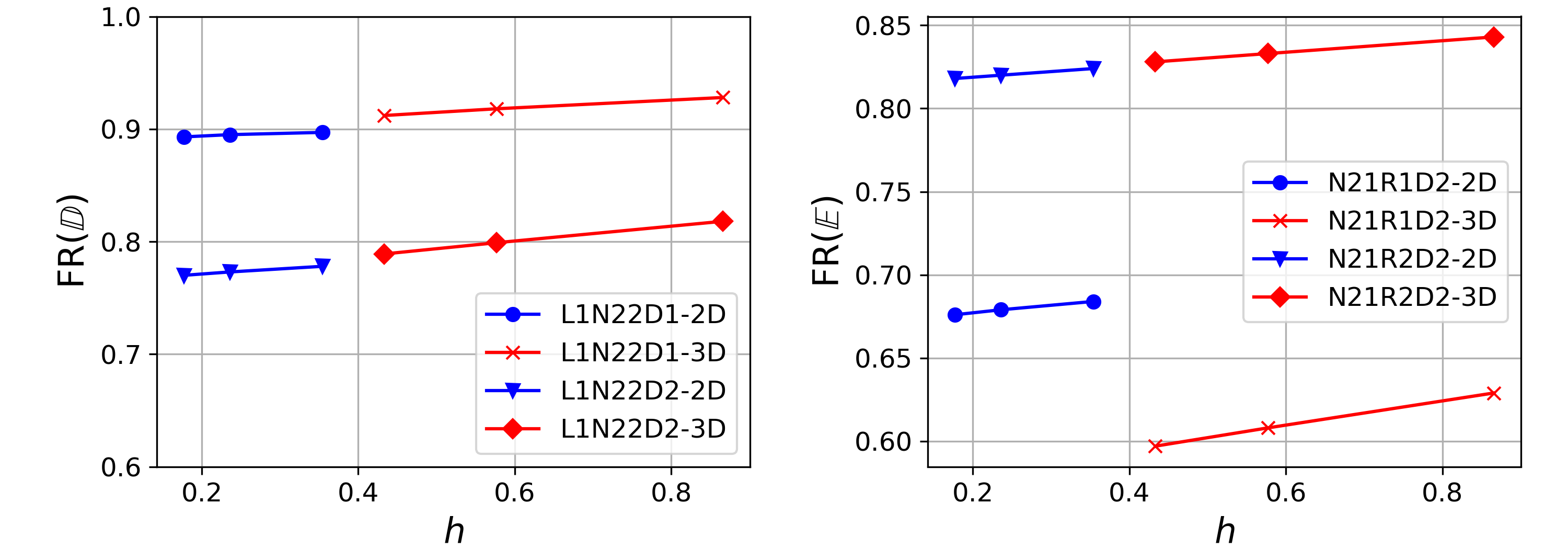}
\end{center}
\vspace*{-0.3in}
\caption{\footnotesize The full-rankness $\mathrm{FR}(\mathbb{M})$ of the matrices of the inf-sup conditions of Section \ref{Sec_StAnalysis} versus the diameter $h$ of the underlying meshes. The full-rankness $\mathrm{FR}(\mathbb{M})$ of a matrix $\mathbb{M}$ is the rank of $\mathbb{M}$ divided by its maximum possible rank. The selected choices of finite elements violate the associated inf-sup condition as the associetd full-rankness is less than $1$. The underlying meshes are shown in Figure \ref{Meshes}.} 
\label{InfSupFR}
\end{figure}

The inf-sup conditions introduced earlier hold only if the corresponding matrices are full rank. Therefore, a simple approach for computationally studying these inf-sup conditions is to compute the rank deficiency of the matrices $\mathbb{S}^{1\mathbf{d}}$, $\mathbb{S}^{\mathrm{D}\mathbf{c}}$, $\mathbb{B}$, $\mathbb{C}$, $\mathbb{D}$, $\mathbb{E}$. Some choices of finite elements that violate these inf-sup conditions are shown in Figure \ref{InfSupFR}. In this figure, the full-rankness of the above matrices are plotted for several choices of finite elements by using $2$D and $3$D meshes of Figure \ref{Meshes}, where the full-rankness $\mathrm{FR}(\mathbb{M})$ of a matrix $\mathbb{M}$ is the rank of $\mathbb{M}$ divided by its maximum possible rank. Thus, $\mathrm{FR}(\mathbb{M})=1$, if $\mathbb{M}$ is full rank and $\mathrm{FR}(\mathbb{M})<1$, otherwise.

The results of Figure \ref{InfSupFR} are computed using the incompressible Neo-Hookean constitutive equation \eqref{IncompNeoHook} with $\mu=1$ near the reference configuration, that is, the matrices are associated to the first Newton's iteration \eqref{NewtIt_IncompFEMs} starting at the reference configuration. The rank deficiency of some finite elements such as $\mathrm{L2R1}$ for $\mathbb{S}^{1\mathbf{d}}$, $\mathrm{N11D2}$ for $\mathbb{S}^{\mathrm{D}\mathbf{c}}$, $\mathrm{L2N11R2D0}$ for $\mathbb{C}$, and $\mathrm{L1N22D1}$ for $\mathbb{D}$, is a simple consequence of the number of columns of the associated submatrix being smaller than the number of its rows as mentioned in Theorem \ref{NecInequal}. For some finite elements such as $\mathrm{L2N11R2D0}$ for $\mathbb{C}$, the associated submatrix is nearly square and as shown in Figure \ref{InfSupFR}, it is also nearly full rank. In this situation, the so-called locking phenomena may occur in the sense that the associated submatrix may represent an injective operator, and thus, it only admits a unique solution.     

Our computations suggest that out of $96$ possible choices of simplicial elements of Figure \ref{FEO012} for displacement, displacement gradient, stress, and pressure, $28$ choices in $2$D and $6$ choices in $3$D satisfy all inf-sup conditions. More specifically, the choices that violate each inf-sup condition are as follow:    

\begin{itemize}
\item[$\bullet$] In $2$D and $3$D, the inf-sup condition \eqref{infsup_S1d} does not hold for the choices of elements $\mathrm{L2R1}$ and $\mathrm{L2B1}$ for displacement and stress.
\item[$\bullet$] The inf-sup condition \eqref{infsup_SDc} does not hold with the following choices of elements for displacement gradient and pressure:
\begin{itemize}
\item[$-$] 2D: $\mathrm{N11D2}$ and $\mathrm{N21D2}$; 
\item[$-$] 3D: $\mathrm{N}i\mathrm{1D1}$ and $\mathrm{N}ij\mathrm{D2}$, $i,j=1,2$.
\end{itemize}
\item[$\bullet$] The inf-sup condition \eqref{infsup_B} does not hold with the following choices of elements for displacement, displacement gradient, and stress:   
\begin{itemize}
\item[$-$] 2D: $\mathrm{L2N11R2}$, $\mathrm{L2N11B2}$, $\mathrm{L2N}ij\mathrm{R1}$, and $\mathrm{L2N}ij\mathrm{B1}$, $i,j=1,2$; 
\item[$-$] 3D: $\mathrm{L2N}ij\mathrm{R1}$ and $\mathrm{L2N}ij\mathrm{B1}$, $i,j=1,2$.
\end{itemize}
\item[$\bullet$] In $2$D and $3$D, the inf-sup condition \eqref{infsup_C} does not hold with the following choices of elements for displacement, displacement gradient, stress, and pressure: $\mathrm{L2N11R2D0}$, $\mathrm{L2N11B2D0}$, $\mathrm{L2N}ij\mathrm{R1D}\ell$, and $\mathrm{L2N}ij\mathrm{B1D}\ell$, $i,j=1,2$, $\ell = 0,1,2$. 
\item[$\bullet$] The inf-sup condition \eqref{infsup_D} does not hold with the following choices of elements for displacement, displacement gradient, and pressure:    
\begin{itemize}
\item[$-$] 2D: $\mathrm{L1N}ij\mathrm{D}k$ and $\mathrm{L2N}ij\mathrm{D2}$, $i,j,k=1,2$; 
\item[$-$] 3D: $\mathrm{L1N}ij\mathrm{D}\ell$ and $\mathrm{L2N}ij\mathrm{D}k$, $i,j,k=1,2$, $\ell = 0,1,2$.
\end{itemize}
\item[$\bullet$] The inf-sup condition \eqref{infsup_E} does not hold with the following choices of elements for displacement gradient, stress, and pressure: 
\begin{itemize}
\item[$-$] 2D: $\mathrm{N}i\mathrm{1R}j\mathrm{D2}$ and $\mathrm{N}i\mathrm{1B}j\mathrm{D2}$, $i,j=1,2$; 
\item[$-$] 3D: $\mathrm{N}i\mathrm{1R}j\mathrm{D}k$, $\mathrm{N}i\mathrm{1B}j\mathrm{D}k$, $\mathrm{N}i\mathrm{2R}j\mathrm{D2}$, and $\mathrm{N}i\mathrm{2B}j\mathrm{D2}$, $i,j,k=1,2$.
\end{itemize}
\end{itemize}

The inf-sup conditions of Section \ref{Sec_StAnalysis} are not independent. For example, the inf-sup conditions \eqref{infsup_S1d} and \eqref{infsup_B} are closely related. Such dependencies can be read off from the above numerical results. Common strategies may be employed to stabilize choices that violate the above inf-sup conditions such as using different underlying meshes for different unknowns, enriching trial and test spaces using bubble functions, and using perturbed formulations \citep{BoBrFo2013}.

\begin{table}[!thb]
\tabcolsep=3.0pt
\caption{Convergence rates $r$ and $L^{2}$-errors of the unit square example: DoF is the number of total degrees of freedom and $(E_{\boldsymbol{U}}, E_{\boldsymbol{K}}, E_{\boldsymbol{P}}, E_{p})=(\|\boldsymbol{U}_{h}-\boldsymbol{U}_{e}\|, \|\boldsymbol{K}_{h}-\boldsymbol{K}_{e}\|, \|\boldsymbol{P}_{h}-\boldsymbol{P}_{e}\|,\|p_{h}-p_{e}\|)$ are the $L^{2}$-errors of the approximate solution ($\boldsymbol{U}_{h}$, $\boldsymbol{F}_{h}$, $\boldsymbol{P}_{h}$, $p_{h}$) with respect to the exact solution \eqref{2Dplate_Exact}. The underlying meshes are the $2$D meshes of Figure \ref{Meshes}.}
\resizebox{0.85\textwidth}{!}{\begin{minipage}{\textwidth}
\centering
\begin{tabular}{cl|cc|cc|cc|cc||cl|cc|cc|cc|cc}
\toprule
FEM & DoF &$E_{\boldsymbol{U}}$ &  & $E_{\boldsymbol{K}}$ &  & $E_{\boldsymbol{P}}$ & & $E_{p}$ &  & FEM & DoF &$E_{\boldsymbol{U}}$ &  & $E_{\boldsymbol{K}}$ &  & $E_{\boldsymbol{P}}$ & & $E_{\boldsymbol{P}}$ &  \\
\midrule
\multirow{4}{*}{\rotatebox[origin=c]{90}{\parbox[c]{2.0cm}{\centering $\mathrm{L1N11R2D0}$}}}
& 546   & 3.50e-3 & \multirow{4}{*}{\rotatebox[origin=c]{90}{\parbox[c]{1.5cm}{\centering $r=2.2$}}}  & 6.27e-2 & \multirow{4}{*}{\rotatebox[origin=c]{90}{\parbox[c]{1.5cm}{\centering $r=1.0$}}}  & 1.45e-1 & \multirow{4}{*}{\rotatebox[origin=c]{90}{\parbox[c]{1.5cm}{\centering $r=1.0$}}} & 8.42e-2 & \multirow{4}{*}{\rotatebox[origin=c]{90}{\parbox[c]{1.5cm}{\centering $r=1.0$}}} &
\multirow{4}{*}{\rotatebox[origin=c]{90}{\parbox[c]{2.0cm}{\centering $\mathrm{L1N11B2D0}$}}} & 722   & 3.64e-3 & \multirow{4}{*}{\rotatebox[origin=c]{90}{\parbox[c]{1.5cm}{\centering $r=2.2$}}}  & 6.25e-2 & \multirow{4}{*}{\rotatebox[origin=c]{90}{\parbox[c]{1.5cm}{\centering $r=1.0$}}}  & 1.45e-1 &\multirow{4}{*}{\rotatebox[origin=c]{90}{\parbox[c]{2.0cm}{\centering $r=1.0$}}} & 8.23e-2 & \multirow{4}{*}{\rotatebox[origin=c]{90}{\parbox[c]{1.5cm}{\centering $r=1.0$}}} \\[5pt]
& 1178   & 1.46e-3 &  & 4.21e-2 &  & 9.61e-2 & & 5.64e-2 & & & 1562   & 1.52e-3 &  & 4.20e-2 &  & 9.59e-2 & & 5.51e-2 &\\[5pt]
& 2050  & 7.86e-4 &  & 3.16e-2 &  & 7.16e-2 & & 4.21e-2 & & & 2722   & 8.18e-4 &  & 3.15e-2 &  & 7.14e-2 & & 4.13e-2 &\\[5pt]
\bottomrule 
\multirow{4}{*}{\rotatebox[origin=c]{90}{\parbox[c]{2.0cm}{\centering $\mathrm{L1N12R2D0}$}}} 
& 786   & 3.50e-3 & \multirow{4}{*}{\rotatebox[origin=c]{90}{\parbox[c]{1.5cm}{\centering $r=2.2$}}}  & 6.27e-2 & \multirow{4}{*}{\rotatebox[origin=c]{90}{\parbox[c]{1.5cm}{\centering $r=1.0$}}}  & 1.45e-1 & \multirow{4}{*}{\rotatebox[origin=c]{90}{\parbox[c]{1.5cm}{\centering $r=1.0$}}} & 8.42e-2 & \multirow{4}{*}{\rotatebox[origin=c]{90}{\parbox[c]{1.5cm}{\centering $r=1.0$}}} &
\multirow{4}{*}{\rotatebox[origin=c]{90}{\parbox[c]{2.0cm}{\centering $\mathrm{L1N12B2D0}$}}} & 962   & 3.64e-3 & \multirow{4}{*}{\rotatebox[origin=c]{90}{\parbox[c]{1.5cm}{\centering $r=2.2$}}}  & 6.25e-2 & \multirow{4}{*}{\rotatebox[origin=c]{90}{\parbox[c]{1.5cm}{\centering $r=1.0$}}}  & 1.45e-1 &\multirow{4}{*}{\rotatebox[origin=c]{90}{\parbox[c]{2.0cm}{\centering $r=1.0$}}} & 8.23e-2 & \multirow{4}{*}{\rotatebox[origin=c]{90}{\parbox[c]{1.5cm}{\centering $r=1.0$}}} \\[5pt]
& 1706  & 1.46e-3 &  & 4.21e-2 &  & 9.61e-2 & & 5.64e-2 & & & 2090   & 1.52e-3 &  & 4.20e-2 &  & 9.59e-2 & & 5.51e-2 &\\[5pt]
& 2978  & 7.86e-4 &  & 3.16e-2 &  & 7.16e-2 & & 4.21e-2 & & & 3650   & 8.18e-4 &  & 3.15e-2 &  & 7.14e-2 & & 4.13e-2 &\\[6pt]
\bottomrule
\multirow{4}{*}{\rotatebox[origin=c]{90}{\parbox[c]{2.0cm}{\centering $\mathrm{L1N21R2D0}$}}}
& 658   & 3.50e-3 & \multirow{4}{*}{\rotatebox[origin=c]{90}{\parbox[c]{1.5cm}{\centering $r=2.2$}}}  & 6.27e-2 & \multirow{4}{*}{\rotatebox[origin=c]{90}{\parbox[c]{1.5cm}{\centering $r=1.0$}}}  & 1.45e-1 & \multirow{4}{*}{\rotatebox[origin=c]{90}{\parbox[c]{1.5cm}{\centering $r=1.0$}}} & 8.42e-2 & \multirow{4}{*}{\rotatebox[origin=c]{90}{\parbox[c]{1.5cm}{\centering $r=1.0$}}} &
\multirow{4}{*}{\rotatebox[origin=c]{90}{\parbox[c]{2.0cm}{\centering $\mathrm{L1N21B2D0}$}}} & 834   & 3.64e-3 & \multirow{4}{*}{\rotatebox[origin=c]{90}{\parbox[c]{1.5cm}{\centering $r=2.2$}}}  & 6.25e-2 & \multirow{4}{*}{\rotatebox[origin=c]{90}{\parbox[c]{1.5cm}{\centering $r=1.0$}}}  & 1.45e-1 &\multirow{4}{*}{\rotatebox[origin=c]{90}{\parbox[c]{2.0cm}{\centering $r=1.0$}}} & 8.23e-2 & \multirow{4}{*}{\rotatebox[origin=c]{90}{\parbox[c]{1.5cm}{\centering $r=1.0$}}} \\[5pt]
& 1418   & 1.46e-3 &  & 4.21e-2 &  & 9.61e-2 & & 5.64e-2 & & & 1802   & 1.52e-3 &  & 4.20e-2 &  & 9.59e-2 & & 5.51e-2 &\\[5pt]
& 2466  & 7.86e-4 &  & 3.16e-2 &  & 7.16e-2 & & 4.21e-2 & & & 3138   & 8.18e-4 &  & 3.15e-2 &  & 7.14e-2 & & 4.13e-2 &\\[5pt]
\bottomrule
\multirow{4}{*}{\rotatebox[origin=c]{90}{\parbox[c]{2.0cm}{\centering $\mathrm{L1N22R2D0}$}}}
& 962   & 3.50e-3 & \multirow{4}{*}{\rotatebox[origin=c]{90}{\parbox[c]{1.5cm}{\centering $r=2.2$}}}  & 6.27e-2 & \multirow{4}{*}{\rotatebox[origin=c]{90}{\parbox[c]{1.5cm}{\centering $r=1.0$}}}  & 1.45e-1 & \multirow{4}{*}{\rotatebox[origin=c]{90}{\parbox[c]{1.5cm}{\centering $r=1.0$}}} & 8.42e-2 & \multirow{4}{*}{\rotatebox[origin=c]{90}{\parbox[c]{1.5cm}{\centering $r=1.0$}}} &
\multirow{4}{*}{\rotatebox[origin=c]{90}{\parbox[c]{2.0cm}{\centering $\mathrm{L1N22B2D0}$}}} & 1138   & 3.64e-3 & \multirow{4}{*}{\rotatebox[origin=c]{90}{\parbox[c]{1.5cm}{\centering $r=2.2$}}}  & 6.25e-2 & \multirow{4}{*}{\rotatebox[origin=c]{90}{\parbox[c]{1.5cm}{\centering $r=1.0$}}}  & 1.45e-1 &\multirow{4}{*}{\rotatebox[origin=c]{90}{\parbox[c]{2.0cm}{\centering $r=1.0$}}} & 8.23e-2 & \multirow{4}{*}{\rotatebox[origin=c]{90}{\parbox[c]{1.5cm}{\centering $r=1.0$}}} \\[5pt]
& 2090   & 1.46e-3 &  & 4.21e-2 &  & 9.61e-2 & & 5.64e-2 & & & 2474   & 1.52e-3 &  & 4.20e-2 &  & 9.59e-2 & & 5.51e-2 &\\[5pt]
& 3650  & 7.86e-4 &  & 3.16e-2 &  & 7.16e-2 & & 4.21e-2 & & & 4322   & 8.18e-4 &  & 3.15e-2 &  & 7.14e-2 & & 4.13e-2 &\\[5pt]
\bottomrule
\multirow{4}{*}{\rotatebox[origin=c]{90}{\parbox[c]{2.0cm}{\centering $\mathrm{L2N12R2D0}$}}}
& 898   & 4.37e-3 & \multirow{4}{*}{\rotatebox[origin=c]{90}{\parbox[c]{1.5cm}{\centering $r=1.9$}}}  & 6.26e-2 & \multirow{4}{*}{\rotatebox[origin=c]{90}{\parbox[c]{1.5cm}{\centering $r=0.8$}}}  & 7.99e-2 & \multirow{4}{*}{\rotatebox[origin=c]{90}{\parbox[c]{1.5cm}{\centering $r=0.9$}}} & 6.83e-2 & \multirow{4}{*}{\rotatebox[origin=c]{90}{\parbox[c]{1.5cm}{\centering $r=1.0$}}} &
\multirow{4}{*}{\rotatebox[origin=c]{90}{\parbox[c]{2.0cm}{\centering $\mathrm{L2N12B2D0}$}}} & 1074   & 4.15e-3 & \multirow{4}{*}{\rotatebox[origin=c]{90}{\parbox[c]{1.5cm}{\centering $r=2.0$}}}  & 5.03e-2 & \multirow{4}{*}{\rotatebox[origin=c]{90}{\parbox[c]{1.5cm}{\centering $r=1.0$}}}  & 7.42e-2 &\multirow{4}{*}{\rotatebox[origin=c]{90}{\parbox[c]{2.0cm}{\centering $r=1.0$}}} & 6.66e-2 & \multirow{4}{*}{\rotatebox[origin=c]{90}{\parbox[c]{1.5cm}{\centering $r=1.0$}}} \\[5pt]
& 1946   & 1.85e-3 &  & 3.43e-2 &  & 4.98e-2 & & 4.45e-2 & & & 2330   & 1.85e-3 &  & 3.36e-2 &  & 4.95e-2 & & 4.44e-2 &\\[5pt]
& 3394  & 1.17e-3 &  & 3.71e-2 &  & 4.33e-2 & & 3.53e-2 & & & 4066   & 1.04e-3 &  & 2.52e-2 &  & 3.71e-2 & & 3.33e-2 &\\[5pt]
\bottomrule
\multirow{4}{*}{\rotatebox[origin=c]{90}{\parbox[c]{2.0cm}{\centering $\mathrm{L2N21R2D0}$}}}
& 770   & 4.37e-3 & \multirow{4}{*}{\rotatebox[origin=c]{90}{\parbox[c]{1.5cm}{\centering $r=1.9$}}}  & 6.26e-2 & \multirow{4}{*}{\rotatebox[origin=c]{90}{\parbox[c]{1.5cm}{\centering $r=0.8$}}}  & 7.99e-2 & \multirow{4}{*}{\rotatebox[origin=c]{90}{\parbox[c]{1.5cm}{\centering $r=0.9$}}} & 6.83e-2 & \multirow{4}{*}{\rotatebox[origin=c]{90}{\parbox[c]{1.5cm}{\centering $r=1.0$}}} &
\multirow{4}{*}{\rotatebox[origin=c]{90}{\parbox[c]{2.0cm}{\centering $\mathrm{L2N21B2D0}$}}} & 946   & 4.15e-3 & \multirow{4}{*}{\rotatebox[origin=c]{90}{\parbox[c]{1.5cm}{\centering $r=2.0$}}}  & 5.03e-2 & \multirow{4}{*}{\rotatebox[origin=c]{90}{\parbox[c]{1.5cm}{\centering $r=1.0$}}}  & 7.42e-2 &\multirow{4}{*}{\rotatebox[origin=c]{90}{\parbox[c]{2.0cm}{\centering $r=1.0$}}} & 6.66e-2 & \multirow{4}{*}{\rotatebox[origin=c]{90}{\parbox[c]{1.5cm}{\centering $r=1.0$}}} \\[5pt]
& 1658   & 1.85e-3 &  & 3.43e-2 &  & 4.98e-2 & & 4.45e-2 & & & 2042   & 1.85e-3 &  & 3.36e-2 &  & 4.95e-2 & & 4.44e-2 &\\[5pt]
& 2882  & 1.17e-3 &  & 3.71e-2 &  & 4.33e-2 & & 3.53e-2 & & & 3554   & 1.04e-3 &  & 2.52e-2 &  & 3.71e-2 & & 3.33e-2 &\\[5pt]
\bottomrule
\multirow{4}{*}{\rotatebox[origin=c]{90}{\parbox[c]{2.0cm}{\centering $\mathrm{L2N22R2D0}$}}}
& 1074   & 4.37e-3 & \multirow{4}{*}{\rotatebox[origin=c]{90}{\parbox[c]{1.5cm}{\centering $r=1.9$}}}  & 6.26e-2 & \multirow{4}{*}{\rotatebox[origin=c]{90}{\parbox[c]{1.5cm}{\centering $r=0.8$}}}  & 7.99e-2 & \multirow{4}{*}{\rotatebox[origin=c]{90}{\parbox[c]{1.5cm}{\centering $r=0.9$}}} & 6.83e-2 & \multirow{4}{*}{\rotatebox[origin=c]{90}{\parbox[c]{1.5cm}{\centering $r=1.0$}}} &
\multirow{4}{*}{\rotatebox[origin=c]{90}{\parbox[c]{2.0cm}{\centering $\mathrm{L2N22B2D0}$}}} & 1250   & 4.15e-3 & \multirow{4}{*}{\rotatebox[origin=c]{90}{\parbox[c]{1.5cm}{\centering $r=2.0$}}}  & 5.03e-2 & \multirow{4}{*}{\rotatebox[origin=c]{90}{\parbox[c]{1.5cm}{\centering $r=1.0$}}}  & 7.42e-2 &\multirow{4}{*}{\rotatebox[origin=c]{90}{\parbox[c]{2.0cm}{\centering $r=1.0$}}} & 6.66e-2 & \multirow{4}{*}{\rotatebox[origin=c]{90}{\parbox[c]{1.5cm}{\centering $r=1.0$}}} \\[5pt]
& 2330   & 1.85e-3 &  & 3.43e-2 &  & 4.98e-2 & & 4.45e-2 & & & 2714   & 1.85e-3 &  & 3.36e-2 &  & 4.95e-2 & & 4.44e-2 &\\[5pt]
& 4066  & 1.17e-3 &  & 3.71e-2 &  & 4.33e-2 & & 3.53e-2 & & & 4738   & 1.04e-3 &  & 2.52e-2 &  & 3.71e-2 & & 3.33e-2 &\\[5pt]
\bottomrule
\end{tabular} 
\label{SquareConvError}
\end{minipage} }
\end{table}

\begin{table}[!thb]
\tabcolsep=3.0pt
\caption{Convergence rates $r$ and $L^{2}$-errors of the unit cube example: DoF is the number of total degrees of freedom and $(E_{\boldsymbol{U}}, E_{\boldsymbol{K}}, E_{\boldsymbol{P}}, E_{p})=(\|\boldsymbol{U}_{h}-\boldsymbol{U}_{e}\|, \|\boldsymbol{K}_{h}-\boldsymbol{K}_{e}\|, \|\boldsymbol{P}_{h}-\boldsymbol{P}_{e}\|,\|p_{h}-p_{e}\|)$ are the $L^{2}$-errors of the approximate solution ($\boldsymbol{U}_{h}$, $\boldsymbol{F}_{h}$, $\boldsymbol{P}_{h}$, $p_{h}$) with respect to the exact solution \eqref{cube_Exact}. The underlying meshes are the $3$D meshes of Figure \ref{Meshes}.}
\resizebox{0.85\textwidth}{!}{\begin{minipage}{\textwidth}
\centering
\begin{tabular}{cl|cc|cc|cc|cc||cl|cc|cc|cc|cc}
\toprule
FEM & DoF &$E_{\boldsymbol{U}}$ &  & $E_{\boldsymbol{K}}$ &  & $E_{\boldsymbol{P}}$ & & $E_{p}$ &  & FEM & DoF &$E_{\boldsymbol{U}}$ &  & $E_{\boldsymbol{K}}$ &  & $E_{\boldsymbol{P}}$ & & $E_{\boldsymbol{P}}$ &  \\
\midrule
\multirow{4}{*}{\rotatebox[origin=c]{90}{\parbox[c]{2.0cm}{\centering $\mathrm{L2N12R2D0}$}}}
& 3243   & 7.35e-3 & \multirow{4}{*}{\rotatebox[origin=c]{90}{\parbox[c]{1.5cm}{\centering $r=2.1$}}}  & 5.57e-2 & \multirow{4}{*}{\rotatebox[origin=c]{90}{\parbox[c]{1.5cm}{\centering $r=1.2$}}}  & 8.21e-2 & \multirow{4}{*}{\rotatebox[origin=c]{90}{\parbox[c]{1.5cm}{\centering $r=1.0$}}} & 5.74e-2 & \multirow{4}{*}{\rotatebox[origin=c]{90}{\parbox[c]{1.5cm}{\centering $r=1.0$}}} &
\multirow{4}{*}{\rotatebox[origin=c]{90}{\parbox[c]{2.0cm}{\centering $\mathrm{L2N12B2D0}$}}} & 4755   & 7.12e-3 & \multirow{4}{*}{\rotatebox[origin=c]{90}{\parbox[c]{1.5cm}{\centering $r=2.0$}}}  & 5.02e-2 & \multirow{4}{*}{\rotatebox[origin=c]{90}{\parbox[c]{1.5cm}{\centering $r=1.1$}}}  & 8.27e-2 &\multirow{4}{*}{\rotatebox[origin=c]{90}{\parbox[c]{2.0cm}{\centering $r=1.0$}}} & 5.74e-2 & \multirow{4}{*}{\rotatebox[origin=c]{90}{\parbox[c]{1.5cm}{\centering $r=1.0$}}} \\[5pt]
& 9993  &  3.15e-3  & & 3.42e-2 &  & 5.51e-2 & & 3.82e-2 & & & 14853   & 3.10e-3 &  & 3.25e-2 &  & 5.53e-2 & & 3.83e-2 &\\[5pt]
& 22611  & 1.75e-3 & & 2.49e-2 &  & 4.13e-2 & & 2.86e-2 & & & 33843   & 1.74e-3 &  & 2.41e-2 &  & 4.15e-2 & & 2.87e-2 &\\[5pt]
\bottomrule
\multirow{4}{*}{\rotatebox[origin=c]{90}{\parbox[c]{2.0cm}{\centering $\mathrm{L2N21R2D0}$}}}
& 2523   & 7.35e-3 & \multirow{4}{*}{\rotatebox[origin=c]{90}{\parbox[c]{1.5cm}{\centering $r=2.1$}}}  & 5.57e-2 & \multirow{4}{*}{\rotatebox[origin=c]{90}{\parbox[c]{1.5cm}{\centering $r=1.2$}}}  & 8.21e-2 & \multirow{4}{*}{\rotatebox[origin=c]{90}{\parbox[c]{1.5cm}{\centering $r=1.0$}}} & 5.74e-2 & \multirow{4}{*}{\rotatebox[origin=c]{90}{\parbox[c]{1.5cm}{\centering $r=1.0$}}} &
\multirow{4}{*}{\rotatebox[origin=c]{90}{\parbox[c]{2.0cm}{\centering $\mathrm{L2N21B2D0}$}}} & 4035   & 7.12e-3 & \multirow{4}{*}{\rotatebox[origin=c]{90}{\parbox[c]{1.5cm}{\centering $r=2.0$}}}  & 5.02e-2 & \multirow{4}{*}{\rotatebox[origin=c]{90}{\parbox[c]{1.5cm}{\centering $r=1.1$}}}  & 8.27e-2 &\multirow{4}{*}{\rotatebox[origin=c]{90}{\parbox[c]{2.0cm}{\centering $r=1.0$}}} & 5.74e-2 & \multirow{4}{*}{\rotatebox[origin=c]{90}{\parbox[c]{1.5cm}{\centering $r=1.0$}}} \\[5pt]
& 7725   & 3.15e-3 &  & 3.42e-2 &  & 5.51e-2 & & 3.82e-2 & & & 12585   & 3.10e-3 &  & 3.25e-2 &  & 5.53e-2 & & 3.83e-2 &\\[5pt]
& 17427  & 1.75e-3 &  & 2.49e-2 &  & 4.13e-2 & & 2.86e-2 & & & 28659   & 1.74e-3 &  & 2.41e-2 &  & 4.15e-2 & & 2.87e-2 &\\[5pt]
\bottomrule
\multirow{4}{*}{\rotatebox[origin=c]{90}{\parbox[c]{2.0cm}{\centering $\mathrm{L2N22R2D0}$}}}
& 3897   & 7.35e-3 & \multirow{4}{*}{\rotatebox[origin=c]{90}{\parbox[c]{1.5cm}{\centering $r=2.1$}}}  & 5.57e-2 & \multirow{4}{*}{\rotatebox[origin=c]{90}{\parbox[c]{1.5cm}{\centering $r=1.2$}}}  & 8.21e-2 & \multirow{4}{*}{\rotatebox[origin=c]{90}{\parbox[c]{1.5cm}{\centering $r=1.0$}}} & 5.74e-2 & \multirow{4}{*}{\rotatebox[origin=c]{90}{\parbox[c]{1.5cm}{\centering $r=1.0$}}} &
\multirow{4}{*}{\rotatebox[origin=c]{90}{\parbox[c]{2.0cm}{\centering $\mathrm{L2N22B2D0}$}}} & 5409   & 7.12e-3 & \multirow{4}{*}{\rotatebox[origin=c]{90}{\parbox[c]{1.5cm}{\centering $r=2.0$}}}  & 5.02e-2 & \multirow{4}{*}{\rotatebox[origin=c]{90}{\parbox[c]{1.5cm}{\centering $r=1.1$}}}  & 8.27e-2 &\multirow{4}{*}{\rotatebox[origin=c]{90}{\parbox[c]{2.0cm}{\centering $r=1.0$}}} & 5.74e-2 & \multirow{4}{*}{\rotatebox[origin=c]{90}{\parbox[c]{1.5cm}{\centering $r=1.0$}}} \\[5pt]
& 11964  & 3.15e-3 &  & 3.42e-2 &  & 5.51e-2 & & 3.82e-2 & & & 16824   & 3.10e-3 &  & 3.25e-2 &  & 5.53e-2 & & 3.83e-2 &\\[5pt]
& 27015  & 1.75e-3 &  & 2.49e-2 &  & 4.13e-2 & & 2.86e-2 & & & 38247   & 1.74e-3 &  & 2.41e-2 &  & 4.15e-2 & & 2.87e-2 &\\[5pt]
\bottomrule
\end{tabular} 
\label{CubeConvError}
\end{minipage} }
\end{table}

\subsection{Convergence Rates}
To study the performance of the mixed finite element method \eqref{IncompFEMs}, we study deformations of a unit square and a unit cube with the incompressible Neo-Hookean constitutive equation \eqref{IncompNeoHook} with $\mu=1$. The underlying meshes are shown in Figure \ref{Meshes}. For the unit square example, we consider the body force and the boundary conditions that induce the displacement and pressure
\begin{equation}\label{2Dplate_Exact}
\boldsymbol{U}_{e}(Y) = \left[\begin{array}{c} \frac{1}{2}Y^{3} + \frac{1}{2}\sin(\frac{\pi}{2}Y) \\ 0 \end{array} \right], \quad p_{e}(Y) = \sin(\frac{\pi}{2}Y).
\end{equation}
As mentioned earlier, $28$ choices out of $96$ possible choices of $2$D elements of Figure \ref{FEO012} satisfy the inf-sup conditions. Numerical results suggest that  $14$ out of these $28$ choices yield poor approximations of stress and pressure. These choices are $\mathrm{L1N}ij\mathrm{R1D0}$, $\mathrm{L1N}ij\mathrm{B1D0}$, $\mathrm{L2N12R2D1}$, $\mathrm{L2N12B2D1}$, $\mathrm{L2N2}j\mathrm{R2D1}$, and $\mathrm{L2N2}j\mathrm{B2D1}$, $i,j=1,2$. Errors and convergence rates of the remaining $14$ choices are shown in Table \ref{SquareConvError}. These results imply that the overall errors and convergence rates strongly depend on the choice of element for stress. Moreover, the stress convergence rates are not optimal in general.

For the unit cube example, we consider the solution
\begin{equation}\label{cube_Exact}
\boldsymbol{U}_{e}(Y) = \left[\begin{array}{c} \frac{1}{4}Y^{3} + \frac{1}{4}\sin(\frac{\pi}{2}Y) \\ 0 \\ 0 \end{array} \right],  \quad p_{e}(Y) = \frac{1}{2}\sin(\frac{\pi}{2}Y).
\end{equation}
Errors and convergence rates for $6$ choices that satisfy all inf-sup conditions are shown in Table \ref{CubeConvError}. Similar to the $2$D example, the results suggest that there is a direct relation between the choice of element for stress and the overall performance of the finite element methods.

\section{Concluding Remarks}
We introduced four-field mixed finite element methods for incompressible nonlinear elasticity, where the independent unknowns are displacement, displacement gradient, the first Piola-Kirchhoff stress tensor, and pressure. Based on a theory for the Galerkin approximation of nonlinear problems, some inf-sup conditions are derived that are sufficient for the local stability of the finite element methods. By considering several choices of simplicial finite elements, we numerically studied the stability and the performance of the associated finite element methods. Although the proposed methods are computationally more expensive than the standard two-field methods for incompressible elasticity, they are potentially useful for accurate approximations of strain and stress as they are independently computed in the solution process. 

For unstable choices of finite elements, one may design a suitable stabilization strategy based on the inf-sup conditions that are violated. Such strategies may  include using different underlying meshes for different unknowns, enriching trial and test spaces using bubble functions, and using perturbed formulations \citep{BoBrFo2013}. We did not study the effect of constitutive equations on the stability of finite element methods. Such an analysis requires physically appropriate assumptions on the elasticity tensor $A^{IJRS}$ or equivalently, on the stored energy function in the case of hyperelastic materials. It is well-known that convexity is not a suitable assumption for nonlinear elasticity. More suitable assumptions include the polyconvexity assumption on the stored energy function or the ellipticity assumption on the elasticity tensor \citep[Section 5.10]{Ci1988}.


\bibliographystyle{unsrtnat}
\bibliography{/home/arzhang/Dropbox/LaTexBibliography/biblio.bib}

%
\end{document}